\documentclass[makeidx]{amsart}
\usepackage{amsmath}
\usepackage{amssymb}
\usepackage{epsfig}
\usepackage{hyperref}
\newtheorem{Proposition}{Proposition}[section]
  \newtheorem{Remark}[Proposition]{Remark}
  
  \newtheorem{Lemma}[Proposition]{Lemma}
  
  \newtheorem{Theorem}{Theorem}[section]
 
 \newtheorem{Definition}[Proposition]{Definition}

\newcommand {\z}{{\noindent}}
\def\blackslug{\hbox{\hskip 1pt \vrule width 4pt height 8pt depth 1.5pt
\hskip 1pt}}
\def\qed{\quad\blackslug\lower 8.5pt\null\par}
 
\def\CC{\mathbb{C}}
 \def\RR{\mathbb{R}}
 \def\NN{\mathbb{N}}

\def\Re{\mathrm{Re}}
\def\Im{\mathrm{Im}}

\def\ds{\text{\raisebox{-2pt}{$\stackrel{_{\displaystyle *}}{*}$}}}

\makeindex
\begin{document}

\author{O. Costin$^1$, and S.  Tanveer$^1$ }\title{Borel summability of
    Navier-Stokes equation in $\mathbb{R}^3$ and small time existence}
\gdef\shortauthors{O.  Costin \& S. Tanveer}
\gdef\shorttitle{Navier-Stokes Equation}
\thanks{$1$.  Department of Mathematics, Ohio State University.}  
 \maketitle

\today

\bigskip

\begin{abstract}

We consider the Navier-Stokes initial value problem,
$$v_t - \Delta v = -\mathcal{P} 
\left [ v \cdot \nabla v \right ] + f ~~,~~v(x, 0) = v_0 (x), ~~x \in
\mathbb{R}^3 $$ where $\mathcal{P}$ is the Hodge-Projection to divergence free
vector fields in the assumption that
 $ \| f
\|_{\mu, \beta} <\infty $ and $\| v_0 \|_{\mu+2, \beta}<\infty$ 
for $\beta\ge 0,\mu>3$, where  
$$  \| {\hat f} (k) \|_{\mu, \beta} = \sup_{k \in \mathbb{R}^3} 
~~e^{\beta |k|} (1+|k|)^\mu | {\hat f} (k) |$$
and ${\hat{f}} (k) = \mathcal{F} [f (\cdot)] (k) $
is the Fourier transform in $x$.

By Borel summation methods we show that there exists a classical solution in
the form
$$ v(x, t) = v_0 + \int_0^\infty e^{-p/t} U(x, p) dp $$
$t\in\CC$, $ \Re \frac{1}{t} > \alpha$, and we estimate $\alpha$ in terms of
$\| {\hat v}_0 \|_{\mu+2, \beta}$ and $ \| {\hat f} \|_{\mu, \beta}$. We show
that $\| {\hat v} (\cdot; t) \|_{\mu+2, \beta} < \infty $.
Existence and $t$-analyticity results are
analogous to Sobolev spaces ones. 

An important feature of the present approach is that continuation of $v$
beyond $t=\alpha^{-1}$ becomes a growth rate question of $U(\cdot, p)$ as $ p
\rightarrow \infty$, $U$ being is a known function. For now, our estimate is
likely suboptimal.

A second result is that we show Borel summability of $v$ for $v_0$ and $f$
analytic. In particular, we obtain Gevrey-1 asymptotics results: $ v \sim v_0
+ \sum_{m=1}^\infty v_m t^m $, where $ |v_m | \le m! A_0 B_0^m$, with $A_0$
and $B_0$ are given in terms of to $v_0$ and $f$ and for small $t$, with
$m(t)=\lfloor B_0^{-1}t^{-1}\rfloor$,
$$ \left| v(x, t) - v_0 (x) - \sum_{m=1}^{m(t)} v_m (x) t^m \right|
\le A_0\,\, m(t)^{1/2}\,\,e^{-m(t)} $$
\end{abstract}

\section{Introduction and main results}

We consider the Navier-Stokes (NS)  initial value problem
\begin{equation}
\label{nseq0}
v_t - \Delta v = -\mathcal{P} [ v \cdot \nabla v ] + f(x)~~,
~~v(x, 0) = v_0 (x), \ \ x\in\RR^3,\ \ t\in\RR^+  
\end{equation}
where $v$ is the fluid velocity and $\mathcal{P}
= I -\nabla \Delta^{-1} (\nabla \cdot )$ 
is the Hodge-Projection 
operator to
the space of divergence free vector fields.  We rescale $v,x$ and $t$ so that
the viscosity is 
one. 
The initial condition $v_0$ and the forcing $f(x)$ are
chosen to be divergence free. We assume $f$ to be  
time-independent
for simplicity, but a time dependent $f$ could be treated similarly. Moreocver,
from the analysis presented here, it will be
clear that similar results
can be obtained for the corresponding periodic problem, 
{\it i.e.} $v (., t) \in \mathbb{T}^3$. 

We first write the equation in the Fourier space. We denote by $\mathcal{F}$
or simply $\hat{ }$ the Fourier transform and ${\hat *}$ is the Fourier
convolution. Since $\nabla \cdot v =0$ we get
\begin{equation}
\label{nseq}
{\hat v}_t + |k|^2 {\hat v} = - i k_j P_k \left [ 
{\hat v}_j {\hat *} {\hat v}  \right ] + {\hat f} ~~,~~{\hat v} (k, 0)=
{\hat v}_0,    
\end{equation}
where as usual a repeated index $j$ denotes summation over
$j\,\,(=1,2, 3)$. If $P_k=\mathcal{F}(\mathcal{P})$ we get 
\begin{equation}
\label{8.0}
P_k \equiv \left ( 1 - \frac{k ( k \cdot )}{|k|^2} \right ),   
\end{equation}
\begin{Definition}
\label{defnorm1}
We introduce the norm $ \| \cdot \|_{\mu, \beta}$ by
\begin{equation}
\label{8.1}
\| {\hat v}_0 \|_{\mu,\beta} 
= \sup_{k \in \mathbb{R}^3}
(1+|k|)^\mu e^{\beta |k|}  
|{\hat v}_0 (k) | ~~,~~{\rm where} ~~{\hat v}_0 (k) = \mathcal{F} [v_0 (\cdot)](k),
\end{equation}
\end{Definition}
We  assume 
$ \| {\hat v}_0 \|_{2 + \mu, \beta} < \infty$, 
$ \| {\hat f} \|_{\mu, \beta} < \infty$
for some 
$\beta \ge 0$ and $\mu > 3$. Clearly, if 
$\beta > 0$, then 
$v_0$ and $f$ are analytic in a strip of width at least $\beta$.

There is considerable mathematical literature for Navier-Stokes equation,
starting with Leray's papers in the 1930s \cite{Leray1}, \cite{Leray2},
\cite{Leray3}.  Global existence and uniqueness are known in 2d (see for
instance \cite{Bertozzi} and reference therein).  However, this is not the
case in 3d.  It is not known whether classical solutions exist globally in
time for arbitrary sized smooth or even analytic initial data.  While weak
solutions in the space of distributions are known to exist since Leray, it is
not known if they are unique or not without additional assumptions.  Only
local existence and uniqueness of classical solutions is known, with a time of
existence inversely proportional to a Sobolev norm of $v_0$. There are
sufficient conditions that guarantee existence for all times \cite{Beale},
\cite{ConstFeffer}, but of course it is unknown whether they are satisfied.
The solution, as long as it exists, is known to be analytic in part of the
right half complex $t$-plane \cite{Masuda}, \cite{Iooss}, \cite{FoiasTem1}.
If space-periodic conditions are imposed, for $v_0 \in H^1 (\mathbb{T}^3) $,
and $f$ analytic, then the solution $v$ becomes analytic in space as well
\cite{FoiasTem}, \cite{DoeringTiti}.

The purpose of this paper is twofold.  One is to introduce Borel transform
techniques (the notions are explained in the sequel) in time for nonlinear
evolution PDEs. After Borel transform, NS becomes an integral equation in $p$,
the Borel dual variable of $1/t$ \footnote{If the equation is first order in
  time and order $n>1$ in space, then $p$ is dual to $t^{-1/(n-1)}$.}  The
integral equation is shown to have a unique solution in an exponentially
weighted space, ${L}^1 (dp \,e^{-\alpha p})$ for some $\alpha>0$. An
important advantage of this formulation is that existence in $t$ of the
evolution PDE is transformed into finding the large $p$-asymptotics of a known
solution to an integral equation (finding $\alpha$).  We do not, in this
paper, find an optimal $\alpha$, but only a rough bound which implies
existence for $t<\alpha^{-1}$.

A second purpose is to show Borel summability of the formal power series in
small $t$ of NS, when initial $v_0$ and $f$ are analytic. This corresponds to
$\beta > 0$ in the norm defined in \ref{defnorm1}.  Borel summability implies
in particular that the formal expansion in powers of $t$,
$$ {\tilde v} (x, t) = v_0 (x) + t v_1 (x)+\cdots $$
where $v_j$ can be found algorithmically, is actually Gevrey-1 asymptotic to
$v$. Borel summability also implies that $\| v_m \|_{\infty} \le 
m! A_0
B_0^m $, where $A_0$ and $B_0$ are determined by $v_0$ and $f$.

Borel summability methods have been used by the authors
\cite{CT1} to 
prove complex sectorial
existence of solutions of a rather general class of nonlinear PDEs 
in $\mathbb{C}^d$ for arbitrary $d$. 
This is in some sense a generalization
of the classical Cauchy-Kowalewski theorem to 
PDEs written as systems  that are first order in time
and higher order in  
space\footnote{Also, Cauchy-Kowalewski  
theorem usually requires a local expansion in all 
{\it all} independent variables. Our methods accommodate
series type expansion in just one variable.}

The main results in this paper are given by the following two theorems. The
results in the first theorem are similar to classical ones, with $\| \cdot
\|_{\mu, \beta}$ replacing
Sobolev norms. 

\begin{Theorem}
\label{localexistence}
If $\| {\hat v}_0 \|_{\mu+2, \beta} < \infty$, $\mu > 3$, $\beta \ge 0$,
NS
has a unique solution $v (\cdot, t)$ such that
$ \| {\hat v} (\cdot, t) \|_{\mu, \beta} < \infty$ 
for 
$\Re \frac{1}{t} > \alpha$. Here $\alpha$ depends on ${\hat v}_0$ 
through (\ref{ensure2}). 

Furthermore, 
${\hat v} (\cdot , t)$ is analytic for $\Re ~\frac{1}{t} > \alpha$ and 
$ \| {\hat v} (\cdot, t) \|_{\mu+2, \beta}<\infty $ for   
$t \in [0, \alpha^{-1} ) $. 
If $\beta > 0$, this implies
that $v$ is analytic in $x$ with the same analyticity width as
$v_0$ and $f$.
\end{Theorem} 

\begin{Remark} 
Sobolev space methods
give local existence of solutions in ${H}^m$
for $t \in [0, T)$, where
$T $ is proportional to $ 1/\| v_0 \|_{{H}^m}$.
In particular, for $m > \frac{7}{2}$, 
these solutions are classical
solutions (the second derivatives are continuous). 
The result in Theorem \ref{localexistence} is similar, but
in a different space. 
The existence time,
$t=\alpha^{-1}$,  involves  
$ \| {\hat v}_0 \|_{j+\mu}$ for $j=0,1,2$
(see (\ref{ensure2})). This solution is classical since
$\| {\hat v} (., t) \|_{\mu+2, \beta} < \infty$ for $\mu > 3$
implies $ v (., t) \in {C}^2 (\mathbb{R}^3)$.  
\end{Remark}   

\begin{Remark}
  If $v_0$ has finite suitable Sobolev norms, it was known that $v $ is
  analytic in $t$ in a region in the right half complex $t$ plane. In
  our setting, $v$ is analytic in $\{t: \Re ~\frac{1}{t} > \alpha \}$ 
  if $\| {\hat v}_0
  \|_{\mu+2, \beta} < \infty$.

\end{Remark}

\begin{Remark}

  Previous results \cite{FoiasTem} show that for space-periodic boundary
  conditions, analytic $f$ and $v_0 \in H^{1}$, the solution $v (\cdot, t)$
  becomes analytic in space, with an analyticity strip improving with time for
  small time. Moreover, for $f=0$, a uniform estimate on the analyticity strip
  width for large time exists under the hypothesis that the local dissipation
  $ \nu \| \nabla v (., t)\|^2_{L^2 (\mathbb{T}^3)}$ is
  bounded \cite{DoeringTiti}.  However, we are not aware of similar results in
  $\mathbb{R}^3$, as is the case in this paper.  For $\beta > 0$, our results
  of Theorem \ref{localexistence} show that the analyticity width is preserved
  for $t \in [0, \frac{1}{\alpha} )$.
\end{Remark}

\begin{Theorem}
\label{T2}
For $\beta > 0$ (analytic initial data)
and $\mu > 3$, the solution $v$
is Borel summable in $1/t$, 
{\it i.e.} there exists 
$U (x, p)$, analytic in a neighborhood of $\RR^+$,
exponentially bounded, and analytic in $x$ for $|\Im ~x | < \beta$
so that
\begin{equation*}
v(x,t) = v_0 (x) + \int_0^{\infty} U(x, p) e^{-p/t} dp 
\end{equation*}
Therefore, in particular, 
as $t\to 0$,
$$ v(x, t) \sim v_0 (x) + \sum_{m=1}^\infty  t^m v_m (x) $$
with 
$$ |v_m (x) | \le m! A_0 B_0^m, $$  
where  $A_0$ and $B_0$ depend on $v_0$ and $f$, through (\ref{12.4.0}), (\ref{ab2cond}) and
(\ref{Bcond})
\end{Theorem}

\begin{Remark}  Borel summability and classical Gevrey-asymptotic
results \cite{Balser} imply for small $t$ that
$$ \left| v(x, t) - v_0 (x) - \sum_{m=1}^{m(t)} v_m (x) t^m \right|
\le A_0 \,\, m(t)^{1/2} e^{-m(t)} $$
where $m(t)=\lfloor B_0^{-1}t^{-1}\rfloor$. Our bounds on 
$B_0$ are likely suboptimal.
Formal 
arguments in the recurrence relation 
of $v_{m+1}$ in terms of $v_m$, $v_{m-1}$,...,$v_1$, 
indicate that
$B$ 
only depends on  $\beta$, but not
on  $\| {\hat v_0} \|_{\mu, \beta}$.
\end{Remark}

\begin{Remark}\label{1.6} For $\beta > 0$ the assumption $\mu > 3$ 
is not
restrictive if  $\beta$ is  consistent with  the analyticity strips of
$v_0$ and $f$. This is because
$(1+|k|)^\mu e^{-{\tilde \beta} |k|} $ is bounded in $k$ for
${\tilde \beta} > 0$.
\end{Remark}

\section{Formulation of Navier Stokes equation: Borel transform}
   
We define ${\hat w}$ by
\begin{equation}
\label{uint}
{\hat v} (k, t) = {\hat v}_0 (k) + t {\hat v}_1 (k) + {\hat w} (k, t) ~~~,
\end{equation}
where
\begin{equation}
\label{4}
{\hat v}_1 (k) = \left ( -|k|^2 {\hat v}_0 - i k_j 
P_k \left [ {\hat v}_{0,j}{\hat*}{\hat v}_0 \right ] \right ) + {\hat f} (k)
\end{equation}
From (\ref{nseq}) we get for ${\hat w}$
\begin{multline}
\label{4.1}
{\hat w}_t + |k|^2 {\hat w} =  
-i k_j 
P_k 
\left [ {\hat v}_{0,j}{\hat *}{\hat w} + 
{\hat w}_j{\hat *}{\hat v}_0 + t {\hat v}_{1,j}{\hat *}{\hat w} + 
t {\hat w}_j {\hat *} {\hat v}_1 + {\hat w}_j {\hat *} {\hat w} \right ] \\ 
- 
 t |k|^2 {\hat v}_1  - i k_j t  
P_k
\left [ {\hat v}_{0,j} {\hat *} {\hat v}_1 + 
{\hat v}_{1,j} {\hat *} {\hat v}_0   
+t {\hat v}_{1,j} {\hat *} {\hat v}_1 \right ] 
\end{multline}

We seek  
a solution as a Laplace transform
\begin{equation}
\label{5}
{\hat w} (k, t) = \int_0^\infty {\hat W} (k, p) e^{-p/t} dp 
\end{equation}
with the property 
$\lim_{p \rightarrow 0+} {\hat W} (k, p) = 0$ and $\lim_{p \rightarrow 0^+}
p {\hat W}_p (k, p) = 0$.  
The Borel transform of (\ref{4.1}), which
is the same as the formal inverse-Laplace transform in $1/t$
gives
in the dual variable $p > 0$, 
\begin{multline}
\label{100.3}
p {\hat W}_{pp} + 2 {\hat W}_p + |k|^2 {\hat W}  
+i k_j 
P_k 
\left [ {\hat v}_{0,j}{\hat *}{\hat W} + 
{\hat W}_j{\hat *}{\hat v}_0 + {\hat v}_{1,j}{\hat *}(1*{\hat W}) + 
(1*{\hat W}_j) {\hat *} v_1 \right ] \\ 
+i k_j P_k  
{\hat W}_j \ds {\hat W}  
+ 
|k|^2 {\hat v}_1  + i k_j 
P_k
\left [ {\hat v}_{0,j} {\hat *} {\hat v}_1 + 
{\hat v}_{1,j} {\hat *} {\hat v}_0 
+p {\hat v}_{1,j} {\hat *} {\hat v}_1 \right ] = 0, 
\end{multline}
where $\ds$ denotes Laplace convolution in $p$, 
followed by Fourier convolution in $k$.

Since the equation
$\mathcal{D} y := [p \partial_p^2 + 2 \partial_p + |k|^2 ] y  = 0$
has explicit independent solutions in terms of
Bessel functions, $y = J_1 (z)/z$ and
$y = Y_1 (z)/z $, where $z=2 |k| \sqrt{p}$ which do not vanish at zero, 
we formally obtain from (\ref{100.3}) by inverting $\mathcal{D}$
the Duhamel formulation
\begin{multline}
\label{100.3.1}
{\hat W} (k, p) = \frac{i k_j \pi}{2 |k| \sqrt{p}} 
\int_0^p \mathcal{G} (z, z') {\hat H}^{[j]} (k, p') dp'  
~~,~~{\rm where}~ \\
\mathcal{G} (z, z') =
z' \left ( -J_1 (z) Y_1 (z') + Y_1 (z) J_1 (z') \right ) 
~~, ~~z=2|k| \sqrt{p} ~,~z'= 2|k| \sqrt{p'} ~,~ 
\end{multline}
and    
\begin{multline}
\label{100.3.2}
{\hat H}^{[j]} = 
-P_k 
\left [ {\hat v}_{0,j}{\hat *}{\hat W} + 
{\hat W}_j{\hat *}{\hat v}_0 + {\hat v}_{1,j}{\hat *}(1*{\hat W}) + 
(1*{\hat W}_j) {\hat *} v_1 \right ] \\ 
-P_k \left [ {\hat W}_j \ds W \right ]  
+ i k_j {\hat v}_1 -
P_k
\left [ {\hat v}_{0,j} {\hat *} {\hat v}_1 + 
{\hat v}_{1,j} {\hat *} {\hat v}_0 
+p {\hat v}_{1,j} {\hat *} {\hat v}_1 \right ]  
\end{multline}
\begin{Remark}\label{bdg}
  $|\mathcal{G}(z,z')|$ is bounded for all real nonnegative $z'\le z$. 
This follows
from standard properties of Bessel functions \cite{Abramowitz}. (The
approximate bound is about  $0.6$.)
\end{Remark}
To obtain stronger results with less regularity
of $v_0$, it is convenient
to introduce ${\hat U} (k, p)$ by:
\begin{equation}
\label{100.3.3}
{\hat W} (k, p) = {\hat U} (k, p) - {\hat v}_1 (k)
\end{equation}  
Substituting (\ref{100.3.3}) into (\ref{100.3.2}), we obtain
\begin{equation}
\label{100.3.4}
{\hat H}^{[j]} (k, p) = {\hat G}^{[j]} (k, p) + i k_j {\hat v}_{1}~~,
~{\rm where} 
~~{\hat G}^{[j]} = -P_k 
\left [ {\hat v}_{0,j}{\hat *}{\hat U} 
+ {\hat U}_j{\hat *}{\hat v}_0 + {\hat U}_j \ds {\hat U} \right ]
\end{equation}
We can further simplify the integral
$\int_0^p \mathcal{G} (z, z') {\hat H}^{[j]} (k, p') dp'$ by
noting that the only solution to 
\begin{equation}
\label{100.3.6}
\mathcal{D} y = -|k|^2 {\hat v}_1 ,  
\end{equation}
satisfying $y (k, 0)=0$, as it is easy to check, is
\begin{equation}
\label{100.3.7}
y (k, p) = -{\hat v}_1 (k) \left ( 1 - 2 \frac{J_1 (z)}{z} \right ) 
~~,~~{\rm where} ~z = 2 |k| \sqrt{p},  
\end{equation}
where we used the fact that $J_1 (z)/z $ is a solution to the
associated homogeneous differential equation and that 
$\lim_{z\rightarrow 0} J_1 (z)/z = 1/2$.
On the other hand, inversion of $\mathcal{D}$
with zero boundary condition at $p=0$ involves
the same kernel $\mathcal{G} (z, z')$.
Writing $-|k|^2 {\hat v}_1 
= i k_j [ ik_j {\hat v}_1 ]$, it follows that
\begin{equation}
\label{100.3.8}
y (k, p) = \frac{i k_j \pi}{2 |k| \sqrt{p}} \int_0^p \mathcal{G} (z, z')
\left [ i k_j {\hat v}_1 (k) \right ] dp' 
\end{equation}
Therefore
\begin{equation}
\label{100.3.5}
\frac{i k_j \pi}{2 k\sqrt{p}} \int_0^p \mathcal{G} (z, z') 
\left [ i k_j {\hat v}_1 ({\hat k}) \right ] = {\hat v}_1 (k) \left ( 
2 \frac{J_1 (z)}{z} - 1 \right )    
\end{equation}
From (\ref{100.3.3}),(\ref{100.3.4}) and (\ref{100.3.5}) we get
\begin{multline}
\label{IntUeqn}
{\hat U} (k, p) = \frac{i k_j \pi}{2 |k| \sqrt{p}} 
\int_0^p \mathcal{G} (z, z') {\hat G}^{[j]} (k, p') dp' 
+ 2 {\hat v}_1 \frac{J_1 \left (2 |k| \sqrt{p} \right )}{2 |k| \sqrt{p}} 
=: \mathcal{N} [{\hat U}] (k, p), 
\end{multline}
where ${\hat G}^{[j]} (k, p)$ is given by (\ref{100.3.4}). 

We will show that $\mathcal{N} $ is contractive in a suitable space, and
hence $ {\hat U} = \mathcal{N} [ {\hat U}]$ has a unique solution. The
solution satisfies 
${\hat U} (0, k) = {\hat v}_1 (k)$, 
${\hat U}$ and ${\hat U}_p $ are bounded for 
$p \in \mathbb{R}^+$ and exponentially bounded at $\infty$.
Then,  
${\hat W} (k, p) = {\hat U} (k, p) - {\hat v}_1 (k) $  satisfies
the integral equation (\ref{100.3.1}) and hence 
the differential equation  
(\ref{100.3}) is satisfied, with 
$ \lim_{p \rightarrow 0} p {\hat W}_p (k, p) = 0$,
$\lim_{p \rightarrow 0} {\hat W} (k, p) =0$, and  
${\hat W}$ and ${\hat W}_p  $ are exponentially bounded at $\infty$. 
Thus the Laplace transform
$ {\hat w} (k, t) = \int_0^\infty e^{-p/t} {\hat W} (k, p) dp $ 
will indeed satisfy  (\ref{4.1}) for sufficiently large
$\Re ~\frac{1}{t}$, and because of the continuity of  ${\hat W}$ at
$p=0$ we have
$\lim_{t \rightarrow 0^+} {\hat w} (k, t) =0$. Thus,
\begin{equation} 
\label{laplace0}
{\hat v} (k, t) = {\hat v}_0 + t {\hat v}_1 +  
\int_0^\infty e^{-p/t} {\hat W} (k, p) dp 
= {\hat v}_0 + \int_0^\infty e^{-p/t} {\hat U} (k, p) dp 
\end{equation} 
solves the NS equation (\ref{nseq}) in the Fourier space, with the given
initial condition.  
Furthermore, the sufficiently rapid decay in $k$ of ${\hat U}$
implies that $ v(x, t) = \mathcal{F}^{-1} [{\hat v} (\cdot, t) ] (x)$ is
indeed a classical solution to (\ref{nseq0}).  It is known
(See {\it e.g.} \cite{Temam}) that classical solutions are
unique; thus 
${\hat v}$ is the only solution to (\ref{nseq0}).

\subsection{Existence of a solution to (\ref{IntUeqn})}

First, we prove some preliminary lemmas.
\begin{Lemma}
\label{lem0.1}
If $\| {\hat v} \|_{\mu,\beta} ~{\rm and}~~ 
\| {\hat w} \|_{\mu,\beta} < \infty$, 
then we have 
\begin{equation}
\| {\hat v} {\hat *} {\hat w} 
\|_{\mu,\beta} \le C_0 \| {\hat v} \|_{\mu,\beta} \| {\hat w} \|_{\mu,\beta}
\end{equation}
where ${\hat *}$ denotes Fourier convolution, 
$$ C_0 (\mu) = 2^{\mu+2} 
\int_{k \in \mathbb{R}^3} \frac{1}{(1+|k|)^\mu} dk= 
\frac{32 \pi 2^{\mu} }{(\mu-1)(\mu-2) (\mu-3)} $$
\end{Lemma}
\begin{proof}
From the definition of $ \| \cdot \|_{\mu, \beta}$, we get
\begin{multline*}
| {\hat v}{\hat *} {\hat w} | \le \| {\hat v} \|_{\mu, \beta} \| {\hat w} 
\|_{\mu, \beta} 
\int_{k' \in \mathbb{R}^3} 
\frac{e^{-\beta (|k'|+|k-k'|)} dk'}{(1+|k'|)^\mu (1+ |k-k'|)^\mu} \\
\le \| {\hat v} \|_{\mu, \beta} \| {\hat w} \|_{\mu, \beta} e^{-\beta |k|}
\int_{k' \in \mathbb{R}^3} 
\frac{dk'}{(1+|k'|)^\mu (1+ |k-k'|)^\mu} 
\end{multline*}
For large $|k|$, we break 
the integral range at $|k'| = |k|/2$.
In the inner ball $|k'| < |k|/2$, we have
$$ \frac{1}{(1+|k'|)^\mu (1+ |k-k'| )^\mu} \le    
\frac{1}{(1+|k'|)^\mu (1+ |k|/2 )^\mu} 
\le \frac{2^\mu}{(1+|k|)^\mu (1+|k'|)^\mu} $$
while, in its complement,
$$ \frac{1}{(1+|k'|)^\mu (1+ |k-k'|)^\mu} \le    
\frac{1}{(1+|k|/2)^\mu (1+ |k-k'|)^\mu } 
\le \frac{2^\mu}{(1+|k|)^\mu (1+|k-k'|)^\mu} $$
Using these estimates, we get for
$\mu > 3$,
\begin{equation}
\label{eqlem0.1}
\int_{k' \in \mathbb{R}^3} \frac{dk'}{(1+|k'|)^\mu (1+ |k-k'|)^\mu} 
\le \frac{C_0}{2 (1+|k|)^\mu} 
\end{equation}
\end{proof}
\begin{Lemma}
\label{lem0.2}
$$
\| P_k \left [ {\hat w}_j {\hat *} {\hat v} \right ] \|_{\mu, \beta} \le 
2 C_0 \| {\hat w}_j \|_{\mu, \beta} \|{\hat v} \|_{\mu, \beta}
$$  
\end{Lemma}
\begin{proof}
It is easily seen from the 
representation of $P_k $ in (\ref{8.0}) that
\begin{equation}
\label{Pbound}
| P_k  {\hat g}  (k) | \le 2 |{\hat g} (k)| 
\end{equation}
Therefore, using (\ref{8.1}), 
$$ \| P_k {\hat g} \|_{\mu,\beta} \le 2 \| g \|_{\mu, \beta} $$
Using Lemma \ref{lem0.1}, with $g = w_j v$,
the proof follows.
\end{proof}
\begin{Lemma}
\label{lemC2}
For $C_2 = 2 \pi C_0 \sup_{z \in \mathbb{R}^+, 0 \le z' \le z} 
| \mathcal{G} (z, z') |$ \footnote{Since sup of $|\mathcal{G}| \approx 0.6$,
we get $C_2 \approx 
\frac{32 (1.2) \pi^2 2^{\mu} }{(\mu-1)(\mu-2) (\mu-3)} $},
with $C_0$ as defined in Lemma \ref{lem0.1},
\label{Nact}
\begin{equation}
\label{N1}
\| \mathcal{N} [{\hat U}] (\cdot , p) \|_{\mu, \beta}  \le 
\frac{C_2}{\sqrt{p}} \int_0^p 
\left \{ \| {\hat U} (\cdot , p') \|_{\mu, \beta} * 
\| U (\cdot , p') \|_{\mu, \beta} 
+ \| v_0 \|_{\mu, \beta} \| {\hat U}(\cdot , p') \|_{\mu, \beta} 
\right \} dp'     
+ \| v_1 \|_{\mu, \beta} 
\end{equation}
\begin{multline}
\label{N2}
\| \mathcal{N} [{\hat U}^{[1]}] (\cdot , p) - 
\mathcal{N} [{\hat U}^{[2]}] (\cdot , p)
\|_{\mu, \beta} \\ 
\le 
\frac{C_2}{\sqrt{p}} \int_0^p  
\left \{ \left ( \| {\hat U}^{[1]} (\cdot , p') \|_{\mu, \beta} 
+ \| {\hat U}^{[2]} (\cdot , p') \|_{\mu, \beta} \right ) * 
\| {\hat U}^{[1]} (\cdot , p') - {\hat U}^{[2]} (\cdot , p')\|_{\mu, \beta} 
\right . \\
\left. + \| v_0 \|_{\mu, \beta} \| {\hat U}^{[1]} (\cdot , p')
-{\hat U}^{[2]} ( \cdot , p') 
\|_{\mu, \beta} \right \} dp' 
\end{multline}
\end{Lemma}
\begin{proof}
From \cite{Abramowitz}, 
$|{J_1 (z)}/{z} | \le 1/2 $ for $z \in
\mathbb{R}^+$ and  therefore
$$\left\| 2{\hat v}_1 (k) {  J_1 (z)}/{z} \right
\|_{\mu, \beta} \le \| {\hat v}_1 \|_{\mu, \beta}$$ 
From 
Lemma \ref{lem0.2}, we have
$$ | \mathcal{P}_k \left \{ {\hat U}_j \ds  {\hat U} \right \}
(k,p ) | 
\le 2 C_0 
\| {\hat U} (\cdot, p) \|_{\mu, \beta} * \| {\hat U} (\cdot, p)\|_{\mu, \beta}  
\frac{e^{-\beta |k|}}{(1+|k|)^\mu} 
$$
Applying Lemma \ref{lem0.2}, we get
$$ \Bigg | \mathcal{P}_k \left \{ {\hat v}_{0_j}{\hat *}{\hat U} (\cdot, p) 
+ {\hat U}_j (\cdot, p) {\hat *} {\hat v}_{0} \right \} \Bigg |
\le 4 C_0 
\| {\hat v}_0  \|_{\mu, \beta} \| {\hat U} (\cdot, p)\|_{\mu, \beta}  
\frac{e^{-\beta |k|}}{(1+|k|)^\mu} 
$$
By Remark~\ref{bdg} and the definition of $\mathcal{N}$ in
(\ref{IntUeqn}), it follows that for
$ C_2 \ge 2 \pi C_0 | \mathcal{G} (z, z') | $
(\ref{N1}) holds.

The second part of the lemma follows by noting that
\begin{equation}
\label{Udiff}
{\hat U}_j^{[1]} \ds  {\hat U}^{[1]}
-{\hat U}_j^{[2]} \ds  {\hat U}^{[2]}
= 
{\hat U}^{[1]}_j \ds  \left ({\hat U}^{[1]} - {\hat U}^{[2]} 
\right )
+ \left( {\hat U}^{[1]}_j - {\hat U}^{[2]}_j \right ) 
 \ds {\hat U}^{[2]} 
\end{equation}
Applying Lemma \ref{lem0.2} to (\ref{Udiff}), we obtain 
\begin{multline*}
\Bigg \| \mathcal{P}_k \left \{
{\hat U}_j^{[1]} \ds  {\hat U}^{[1]} (\cdot, p)
-{\hat U}_j^{[2]} \ds  {\hat U}^{[2]} (\cdot, p) \right \}
\Bigg \|_{\mu, \beta} 
\\ \le  2 C_0 
\| {\hat U}^{[1]} (\cdot, p) \|_{\mu, \beta} * \| {\hat U}^{[1]} (\cdot, p) - 
{\hat U}^{[2]} (\cdot, p) 
\|_{\mu, \beta} \\
+  
2 C_0 \| {\hat U}^{[2]} (\cdot, p) \|_{\mu, \beta} * 
\| {\hat U}^{[1]} (\cdot, p) - {\hat U}^{[2]} (\cdot, p) 
\|_{\mu, \beta} ,
\end{multline*}
from which (\ref{N2}) follows easily.
\end{proof}

It is convenient to define a number of
different norms for functions of $(k, p)$ on 
$\mathbb{R}^3 \times \left ( \mathbb{R}^+ \cup \{0\} \right )$
\begin{Definition}
For $\alpha \ge 1$, we define
\begin{equation}
\label{8}
\| {\hat f} \|^{(\alpha)} = \sup_{p \ge 0}
(1+p^2) e^{-\alpha p} | {\hat f} (\cdot, p) |_{\mu, \beta}
\end{equation}
We define $\mathcal{A}^\alpha$ to be the
Banach-space of 
continuous functions of $(k, p)$ for $k \in \mathbb{R}^3$ and
$p \in [0, \infty)$ for which
$ \| . \|^{(\alpha)} < \infty$.
It is also convenient to consider the Banach space $\mathcal{A}_1^{\alpha} $
of
locally integrable ($L^1_{loc} $) functions for $p \in [0, L)$  
on $\mathbb{R}^+$, and continuous in $k \in \mathbb{R}^3$ 
such that
\begin{equation}
\label{8.0.0}
\| {\hat f} \|_1^{(\alpha)} = \int_0^L e^{-\alpha p} 
\| {\hat f} (\cdot, p) \|_{\mu, \beta} dp~~  
< \infty ~,\end{equation}
where $L$ is allowed to be finite or $\infty$. 
It is also convenient to define 
$\mathcal{A}_L^{\infty}$ to be the Banach space of 
continuous functions of $(k, p)$ on
$\mathbb{R}^3 \times [0, L]$
such that 
\begin{equation}
\label{8.0.0.0}
\|{\hat f} \|_{L}^{(\infty)} 
= \sup_{p \in [0, L]} \| {\hat f} (\cdot, p) \|_{\mu, \beta} < \infty
\end{equation}
\end{Definition}

\begin{Lemma}
\label{lemBanach}
For ${\hat f}, {\hat g} \in 
\mathcal{A}^\alpha, \mathcal{A}_1^{\alpha}$ or $\mathcal{A}_L^{\infty}$, 
we have the following
the following Banach algebra properties:
$$ \| {\hat f} \ds  {\hat g} \|^{(\alpha)} 
\le M_0 \| {\hat f} \|^{(\alpha)}
\| {\hat g} \|^{(\alpha)}, ~~{\rm where}~M_0 \approx 3.76 \cdots $$
$$ \| {\hat f} \ds  {\hat g} \|_1^{(\alpha)} 
\le \| {\hat f} \|_1^{(\alpha)} 
\| {\hat g} \|_1^{(\alpha)}, $$   
$$ \| {\hat f} \ds {\hat g} \|_L^{(\infty)} 
\le L \| {\hat f} \|_L^{(\infty)} 
\| {\hat g} \|_L^{(\infty)} $$   
\end{Lemma}
\begin{proof} 
In the following, 
we take $u(p) = \| {\hat f}(\cdot, p) \|_{\mu, \beta}$ and
$v (p) = \| {\hat g} (\cdot, p) \|_{\mu, \beta}$. 
We observe that
\begin{multline*}
\int_0^L u(s) v(p-s) ds \le e^{\alpha p}
\left ( \sup_{p \in \mathbb{R}^+} (1+p^2) e^{-\alpha p} u (p) \right ) 
\left ( \sup_{p \in \mathbb{R}^+} (1+p^2) e^{-\alpha p} v (p) \right ) \times
\\
\int_0^p \frac{ds}{(1+s^2)[1+(p-s)^2]} 
\end{multline*}
The first part of the lemma follows
since  \cite{CT1}
$$\int_0^p \frac{ds}{(1+s^2)[1+(p-s)^2]} \le \frac{M_0}{1+p^2} $$ with
$M_0 = 3.76 \cdots$. 
For the second part note that
\begin{multline}
 \int_0^L e^{-\alpha p} \int_0^p u(s) v (p-s) ds 
\\=
\int_0^L \int_0^p e^{-\alpha s} e^{-\alpha (p-s)} u(s) v (p-s) ds 
\le \int_0^L e^{-\alpha s} u(s) ds  \int_0^L e^{-\alpha \tau} v(\tau)  
\end{multline}
The third part follows from the fact that for $p \in [0, L]$
$$ \int_0^p |u(s)| |v(p-s)| \le 
\left \{ \sup_{p\in [0,L]} | u (p)|  
\right \}
\left ( \sup_{p\in [0,L]} | v (p) | \right ) L
$$\end{proof}

\begin{Lemma} 
\label{NUnormbound}
On $\mathcal{A}_1^{\alpha}$, 
the operator $\mathcal{N} $, defined in (\ref{IntUeqn}), 
satisfies 
the following inequalities, with $C_2$ defined in Lemma 
\ref{lemC2}:
\begin{equation}
  \| \mathcal{N} [{\hat U}] \|_1^{(\alpha)} \le  
C_2 \sqrt{\pi} \alpha^{-1/2} 
\left \{ \left ( \| {\hat U} \|_1^{(\alpha)} \right )^2 + 
\| {\hat v}_0 \|_{\mu, \beta} \| {\hat U} \|_1^{(\alpha)} 
\right \} + \alpha^{-1} 
\| {\hat v}_1 \|_{\mu, \beta}
\end{equation}
\begin{multline}
   \| \mathcal{N} [{\hat U}^{[1]} ] - \mathcal{N}
[ {\hat U}^{[2]}] \|_1^{(\alpha)}\\ \le  
C_2 \sqrt{\pi} \alpha^{-1/2} 
\left \{ \left ( \| {\hat U}^{[1]} \|_1^{(\alpha)} +
\| {\hat U}^{[2]} \|_1^{(\alpha)} 
\right ) \| {\hat U}^{[1]} - {\hat U}^{[2]} \|_1^{(\alpha)} + 
\| {\hat v}_0 \|_{\mu, \beta} \| {\hat U}^{[1]}-{\hat U}^{[2]} \|_1^{(\alpha)} 
\right \}
\end{multline}
while in $\mathcal{A}_L^{\infty}$, we have
\begin{equation}
\label{NULinf}
\| \mathcal{N} [{\hat U}] \|_L^{(\infty)} \le  
C_2 L^{1/2} \left \{ L \left ( \| {\hat U} \|_L^{(\infty)} \right )^2 + 
\| {\hat v}_0 \|_{\mu, \beta} \| {\hat U} \|_L^{(\infty)} \right \} +   
\| {\hat v}_1 \|_{\mu, \beta} 
\end{equation}

\begin{multline}
   \| \mathcal{N} [{\hat U}^{[1]} ] - \mathcal{N}
[ {\hat U}^{[2]}] \|_L^{(\infty)}\\ \le  
C_2 L^{1/2} \left \{ L \left ( \| {\hat U}^{[1]} \|_L^{(\infty)} + 
\| {\hat U}^{[2]} \|_L^{(\infty)} 
\right ) \| {\hat U}^{[1]} - {\hat U}^{[2]} \|_L^{(\infty)} + 
\| {\hat v}_0 \|_{\mu, \beta} \| {\hat U}^{[1]}-{\hat U}^{[2]} 
\|_L^{(\infty)} 
\right \} 
\end{multline}   
\end{Lemma}
\begin{proof}
For the space $\mathcal{A}_1^\alpha$, for any $L > 0$, including
$L = \infty$, we note that
$$ \int_0^L e^{-\alpha p} \| {\hat v}_1 \|_{\mu, \beta} dp \le 
\alpha^{-1} \| {\hat v}_1 \|_{\mu, \beta}, $$
while
$$ \int_0^L p^{-1/2} e^{-\alpha p} dp 
\le \Gamma \left(\frac{1}{2} \right) \alpha^{-1/2}=\sqrt{\pi} \alpha^{-1/2} $$  
Furthermore, we note that for $u (p') \ge 0$ we have 
\begin{multline}
\int_0^L e^{-\alpha p} p^{-1/2} \left ( \int_0^p u(p') dp' \right )    
= \int_0^L u(p') e^{-\alpha p'} \left ( \int_{p'}^L p^{-1/2} 
e^{-\alpha (p-p')} dp
\right ) dp' \\
\le 
\int_0^L e^{-\alpha p'} u(p') \int_0^L s^{-1/2} e^{-\alpha s} ds dp'  
\end{multline}
Therefore, it follows from (\ref{N1}) that 
\begin{equation}
\int_0^L e^{-\alpha p} \| \mathcal{N} [{\hat U} ] (\cdot, p) \|_{\mu, \beta} dp
\le C_2 \sqrt{\pi} \alpha^{-1/2} \left ( \left [ \| {\hat U} \|_1^{(\alpha)} 
\right ]^2   
+ \| v_0 \|_{\mu, \beta} \| {\hat U} \|_{1}^{(\alpha)} \right ) +  
\alpha^{-1} \| v_1 \|_{\mu, \beta}
\end{equation}
Furthermore, from (\ref{N2}), it follows that 
\begin{multline*}
\int_0^L \| \mathcal{N} [{\hat U}^{[1]} ] - \mathcal{N}
[ {\hat U}^{[2]}] \|_{\mu, \beta} e^{-\alpha p} dp\\ \le  
C_2 \sqrt{\pi} \alpha^{-1/2} \left \{ 
\left ( \| {\hat U}^{[1]} \|_1^{(\alpha)}
+
\| {\hat U}^{[2]} \|_1^{(\alpha)} 
\right ) \| {\hat U}^{[1]} - {\hat U}^{[2]} \|_1^{(\alpha)} \right .\\
\left. + 
\| {\hat v}_0 \|_{\mu, \beta} \| {\hat U}^{[1]}-{\hat U}^{[2]} \|_1^{(\alpha)} 
\right \} 
\end{multline*}
Hence the first part of the lemma follows.

\z For the second part, we first note that for any $p \in [0, L]$ we have
\begin{equation}
\label{inequal1}
\left| p^{-1/2} \int_0^p u(p') dp' \right| \le \sup_{p \in [0, L]} |u(p) | \sqrt{L} 
\end{equation}
We note that  
\begin{equation}
\label{inequal2}
\left| \int_0^p y_1 (s) y_2 (p-s) ds \right| \le L 
\left ( \sup_{p \in [0, L]}  |y_1 (p)| \right )
\left ( \sup_{p \in [0, L]}  |y_2 (p)| \right ) 
\end{equation}
Taking
$$u(p) = \| {\hat U} (\cdot, p) \|_{\mu, \beta} *  
\| {\hat U} (\cdot, p) \|_{\mu, \beta} + \| v_0 \|_{\mu, \beta} 
\| {\hat U} (\cdot, p) \|_{\mu, \beta} $$ 
$$y_1 (p) = y_2 (p) = \| {\hat U} (\cdot, p) \|_{\mu, \beta}$$ 
(\ref{NULinf}) 
follows from  (\ref{N1}). 
To bound $\mathcal{N} [{\hat U}^{[1]}] - \mathcal{N} [{\hat U}^{[2]} ]$
in $\mathcal{A}_L^{\infty}$, we take
\begin{multline}
 u(p) = \left ( \| {\hat U}^{[1]} (\cdot, p) \|_{\mu, \beta} +
\| {\hat U}^{[2]} (\cdot, p) \|_{\mu, \beta} \right ) *  
\| {\hat U}^{[1]} (\cdot, p) - {\hat U}^{[2]} (\cdot, p) \|_{\mu, \beta}  
\\+ \| v_0 \|_{\mu, \beta} 
\| {\hat U}^{[1]} (\cdot, p) - {\hat U}^{[2]} (\cdot,p) \|_{\mu, \beta} 
\end{multline}

$$y_1 (p) = 
\left ( \| {\hat U}^{[1]} (\cdot, p) \|_{\mu, \beta} +
\| {\hat U}^{[2]} (\cdot, p) \|_{\mu, \beta} \right )   
;\ y_2 (p) = 
\| {\hat U}^{[1]} (\cdot, p) 
- {\hat U}^{[2]} (\cdot, p) \|_{\mu, \beta} $$ 
in  (\ref{inequal1}) and (\ref{inequal2}). 
 The proof now follows from (\ref{N2}).
\end{proof}
\begin{Lemma}
\label{inteqn}
Equation 
(\ref{IntUeqn}) has a unique solution in $\mathcal{A}_1^{\alpha}$ for 
any
$L > 0$ (including $L=\infty$) 
in a ball of size
$ 2\alpha^{-1} \| {\hat v}_1 \|_{\mu, \beta}$, for $\alpha$   large enough to ensure
\begin{equation}
\label{ensure2}
2 C_2 \sqrt{\pi} \alpha^{-1/2} 
\left ( \| v_0 \|_{\mu, \beta} + 2 \alpha^{-1} \| v_1 \|_{\mu, \beta} 
\right ) < 1,  
\end{equation}
where 
$C_2 \approx
\frac{32 (1.2) \pi^2 2^{\mu} }{(\mu-1)(\mu-2) (\mu-3)} $ is the same
as in Lemma \ref{lemC2}.
Furthermore, this solution belongs to $\mathcal{A}_L^{\infty}$ for $L$
small enough so that
\begin{equation}
\label{ensure3}
2 C_2 L^{1/2}  
\left ( \| v_0 \|_{\mu, \beta} + 2 L  \| v_1 \|_{\mu, \beta} 
\right ) < 1,  
\end{equation}
In particular, $\lim_{p \rightarrow 0} {\hat U} (k, p) = {\hat v}_1 (k)$.
 Also, ${\hat W} (k, p) = {\hat U} (k, p) - {\hat v}_1 (k)$
is the unique solution to (\ref{100.3}) which is zero at $p=0$.
\end{Lemma}
\begin{proof}
The estimates of Lemma 
\ref{NUnormbound}
imply that $\mathcal{N}$ maps a ball of size 
$ 2 \alpha^{-1} \| v_1 \|_{\mu, \beta} $ in $\mathcal{A}_1^{\alpha}$
back to itself and that $\mathcal{N}$ is contractive in that ball
when $\alpha$ satisfies (\ref{ensure2}).
From Lemma \ref{NUnormbound} in
space $\mathcal{A}_L^{\infty}$, it follows that $\mathcal{N} $ maps a ball
of size $2 \| v_1 \|_{\mu, \beta} $ to itself and that $\mathcal{N}$
is also contractive in this ball if $L$ is small
enough to ensure (\ref{ensure3}). Thus, there is a unique solution in
this ball. Since $\mathcal{A}_L^{\infty} 
\subset \mathcal{A}_1^{\alpha} $, it follows that the
solutions are in fact the same.

Using Lemma \ref{NUnormbound}, with 
$ {\hat  U}^{[1]} = {\hat U}$ and ${\hat U}^{[2]} = 0$, we obtain from   
(\ref{IntUeqn}), 
$$ \Big \| {\hat U} (k, p) - {\hat v}_1 (k) \frac{2 J_1 (z)}{z} 
\Big \|_L^{(\infty)}
\le C_2 L^{1/2}  \left ( L \left [ \| {\hat U} \|_L^{(\infty)} \right ]^2 +    
\| {\hat v}_0 \|_{\mu, \beta} \| {\hat U} \|_L^{(\infty)} \right )   
$$ 
Since $\| {\hat U} \|_L^{(\infty)} < 2 \| {\hat v}_1 \|_{\mu, \beta}$, 
it follows
that as $L \rightarrow 0$, 
$$ \| {\hat U} (k, p) - 2{\hat v}_1 (k)  J_1 (z)/{z} 
\|_L^{(\infty)} \rightarrow 0$$
Since
$\lim_{z \rightarrow 0} 
2 J_1 (z)/{z} = 1$, it follows that for fixed $k$,
$\lim_{p \rightarrow 0}  {\hat U} (k, p) = {\hat v}_1 (k)$.
By construction, ${\hat U}$ satisfies (\ref{IntUeqn}) iff
${\hat W} = {\hat U} - {\hat v}_1 $ satisfies 
(\ref{100.3.1}). From the properties of 
$\mathcal{G}$ and ${\hat H}^{[j]}$,
it follows that ${\hat W} $ will indeed satisfy (\ref{100.3})
and that it is the only solution which is zero at $p=0$.
\end{proof} 

\begin{Proposition}
\label{inteqnA}
If $\alpha$ is large enough so that (\ref{ensure2}) holds, then 
for an absolute constant  $C_3 > 0$, the solution
${\hat U} (k, p)$ in Lemma \ref{inteqn} and its $p$-derivative satisfy
$$
|{\hat U} (k, p) | \le 
\frac{2 e^{-\beta |k|+\alpha p} \| {\hat v}_1 \|_{\mu, \beta} }{ 
(1+|k|)^\mu} 
$$
\begin{equation*}
| {\hat U}_p (k, p) | 
\le \frac{C_3 e^{-\beta |k|} \| {\hat v}_1 \|_{\mu, \beta}}{(1+|k|)^\mu} 
\left \{ \frac{\sqrt{\alpha}}{C_2} |k| e^{\alpha p} + |k|^2 \right \} 
\end{equation*}
In particular, ${\hat U} \in \mathcal{A}^{\alpha'}$ for any  
$\alpha' > \alpha$, and 
$$
|{\hat U} (k, p) | \le 
\left ( \sup_{p \in \mathbb{R}^+} (1+p^2) e^{-(\alpha'-\alpha) p} \right ) 
\frac{2 e^{-\beta |k|+\alpha' p} \| {\hat v}_1 \|_{\mu, \beta} }{ 
(1+p^2) (1+|k|)^\mu} 
$$

\end{Proposition}
\begin{proof}
With $L = L_0 =\alpha^{-1}$, then
(\ref{ensure3})  holds, and therefore
${\hat U} \in \mathcal{A}_{L_0}^{\infty}$.
For $p \in [0, L_0]$, we obtain
\begin{equation}
\label{eqinside}
e^{-\alpha p} \| {\hat U} (\cdot, p) \|_{\mu, \beta}
< 2 e^{-\alpha p} \| {\hat v}_1 \|_{\mu, \beta} 
\end{equation}
We now consider $p \in [L_0, \infty)$.
We define
$$ y(p) = \|{\hat U} (\cdot , p) \|_{\mu, \beta} * 
\|{\hat U} (\cdot, p) \|_{\mu, \beta} + \| {\hat v}_0 \|_{\mu, \beta} 
\|{\hat U} (\cdot, p) \|_{\mu, \beta} $$ 
We note that
\begin{equation}
\label{newA}
\left\lvert  \frac{1}{\sqrt{p}} e^{-\alpha p} \int_0^p y (p') dp' \right\rvert   
\le L_0^{-1/2} \left\lvert  \int_0^p e^{-\alpha p'} y(p') dp' \right\rvert \le 
\alpha^{1/2} \| y \|_1^{(\alpha)} 
\end{equation}     
From (\ref{IntUeqn}) and (\ref{ensure2}), 
it follows that for $ p \in [L_0, \infty)$
\begin{multline}\label{ab}
|{\hat U} (k, p) | \le \frac{e^{-\beta |k|+\alpha p} }{
(1+|k|)^\mu} 
\left \{ C_2 \alpha^{1/2} \left ( \| {\hat U} \|_1^{(\alpha)} \right )^2 + 
C_2 \alpha^{1/2} \| {\hat v}_0 \|_{\mu, \beta} \| {\hat U} \|_1^{(\alpha)} 
+ e^{-\alpha p} \| {\hat v}_1 \|_{\mu, \beta} \right \}  \\
\le  
\frac{2 e^{-\beta |k|+\alpha p} }{ 
(1+|k|)^\mu} 
\| {\hat v}_1 \|_{\mu, \beta} 
\end{multline}
By (\ref{eqinside}), (\ref{ab}) holds for $p \in [0, L_0]$ as
well; hence the 
bound for $|{\hat U} |$ follows. 
For $\alpha' > \alpha$, $ \| {\hat U} \|^{(\alpha')} < \infty$  because
$e^{-(\alpha' -\alpha) p} (1+p^2)$ is bounded if $\alpha' > 
\alpha$. 

Since ${\hat U}$ is a solution to 
(\ref{IntUeqn}), 
 differentiation with respect to $p$ implies that 
\begin{multline*}
{\hat U}_p (k, p) = 
{\hat v}_1 (k)
\left ( \frac{J_1 (z)}{z} \right )^\prime \frac{4 |k|^2}{z}
~+~ \frac{i k_j \pi }{p} \int_0^p \left \{ \mathcal{G}_z (z, z') - 
\frac{\mathcal{G} (z, z')}{z} \right \} {\hat G}^{[j]} (k, p') dp'
\end{multline*}
Since the functions
$\mathcal{G}_z (z, z')$, $ \mathcal{G} (z, z')/z $ and 
$ z^{-1} \left ( {J_1 (z)}/{z} \right )^\prime $
are easily checked to be bounded for  $z\ge z' \in \mathbb{R}^+$, 
there exists
$C_3 > 0$, independent of any parameter, so that
\begin{multline*}
|{\hat U}_p (k, p) | \le  \frac{C_3 |k|}{p} \left| \int_0^p
|{\hat G} | (k, p') | dp' + C_3 |k|^2 |{\hat v}_1 (k) \right|   
\le 
\frac{C_3 |k| e^{-\beta |k|}}{
(1+|k|)^\mu} \\ \times \left [ \frac{1}{p} \int_0^p 
\left ( \| U (\cdot, p') \|_{\mu, \beta} * \|U (\cdot, p')\|_{\mu, \beta} +  
\| {\hat v}_1 \|_{\mu, \beta} \| U (\cdot, p') \|_{\mu, \beta} \right ) dp' 
+ |k| \| {\hat v}_1 \|_{\mu, \beta} \right ]   
\end{multline*}
For $ p \in [0, L_0]$, with $L=L_0 = \frac{1}{\alpha}$ 
satisfying 
(\ref{ensure3}), we have
\begin{multline*}
|{\hat U}_p (k, p)| \le 
\frac{C_3 |k| e^{-\beta |k|}}{(1+|k|)^\mu} 
\left [ \left \{ L_0 \left ( \| U \|_{L_0}^{(\infty)} \right )^2 + 
\| {\hat v}_0 |_{\mu, \beta} \| U \|_{L_0}^{(\infty)} \right \}   
+ |k| \| {\hat v}_1 \|_{\mu, \beta} \right ] \\
\le   
\frac{C_3 e^{-\beta |k|}}{(1+|k|)^\mu}  
\left ( \frac{\sqrt{\alpha}}{C_2} |k| + 
|k|^2 \right ) \| {\hat v}_1 \|_{\mu, \beta} 
\end{multline*}
For $ p \in [L_0, \infty)$ and $\alpha$ satisfying (\ref{ensure2}), we have 
\begin{multline*} 
|{\hat U}_p (k, p) | \le 
\frac{C_3 |k| e^{-\beta |k|+\alpha p}}{L_0 (1+|k|)^\mu} 
\left [ \left \{  \left ( \| U \|_{1}^{(\alpha)} \right )^2 + 
\| {\hat v}_0 |_{\mu, \beta} \| U \|_{1}^{(\alpha)} \right \}   
+ |k| L_0 e^{-\alpha p} \| {\hat v}_1 \|_{\mu, \beta} \right ] \\
\le   
\frac{C_3 e^{-\beta |k|} \| {\hat v}_1 \|_{\mu, \beta}}{(1+|k|)^\mu} 
\left \{ \frac{\sqrt{\alpha}}{C_2} |k| e^{\alpha p} + |k|^2 \right \} 
\end{multline*}
Continuity of ${\hat U}$ in $p$ follows from the
boundedness of ${\hat U}_p$ for $p \in \mathbb{R}^+$ for
fixed $k$.
\end{proof}

\begin{Lemma}
\label{localexistence1}
Let $\| {\hat v}_0 \|_{\mu+2, \beta} < \infty$ and $\|{\hat f} \|_{\mu, \beta}
< \infty$, with $\mu > 3$, $\beta \ge 0$. Then NS has a unique solution with $
\| {\hat v} (\cdot, t) \|_{\mu, \beta}<\infty$ and ${\hat v} (\cdot, t)$
analytic in $t$ for $\Re \frac{1}{t} > \alpha$, where $\alpha$ depends on the
initial data (see (\ref{ensure2})).  For $\beta > 0$, this implies $v$ is
analytic in $x$ in the same analyticity strip as $v_0,f$.
\end{Lemma}

\begin{proof}
  From (\ref{4}) we see that $ \|{\hat v}_1 \|_{\mu, \beta} < \infty$, since
\begin{equation}
\| {\hat v}_1 \|_{\mu, \beta} \le  \| {\hat v}_0 \|_{\mu +2, \beta} 
+ 2 C_0 \| {\hat v}_0  \|_{\mu, \beta} \| {\hat v}_0 \|_{\mu+1, \beta} 
+ \| {\hat f} \|_{\mu, \beta}
\end{equation}   
Therefore, when $\alpha$ is large enough to ensure (\ref{ensure2}), it follows
that ${\hat U} (k, \cdot)$ and ${\hat W} (k, .) \equiv {\hat U} (k\cdot, .) -
{\hat v}_1 (k)$ are in ${L}^1(e^{-\alpha p}dp )$.  From Lemma \ref{inteqn},
it follows that $\lim_{p \rightarrow 0} {\hat W} (k, p) =0$ and Proposition
\ref{inteqnA} implies ${\hat W}_p (k, p)$ (same as $ {\hat U}_p (k, p)$) is
bounded for $p \in \mathbb{R}^+$ and hence $\lim_{p \rightarrow 0^+} p {\hat
  W}_p = 0$.  Since ${\hat U}$ satisfies (\ref{IntUeqn}), it follows that
${\hat W}$ will satisfy (\ref{100.3.1}) and hence (\ref{100.3}). 
For 
$\Re\,t^{-1}>\alpha$, we take the Laplace transform of ({\ref{100.3}) in $p$,
  using the fact $ \partial_p [p {\hat W}] $ and $p {\hat W} $ vanish at
  $p=0$. There is no contribution at $\infty$ because of boundedness of $
  e^{-\alpha p } \left ( |{\hat W} | + {\hat W}_p \right )$ which follows from
  Proposition \ref{inteqnA}.  It can be checked that}
  ${\hat w} (k, t) = \int_0^\infty
  {\hat W} (k, p) e^{-p/t} dp$ satisfies (\ref{4.1}).
  Therefore,
$$ {\hat v} (k, t) = {\hat v}_0 + t {\hat v}_1 + 
\int_0^\infty {\hat W} (k, p) e^{-p/t} dp = {\hat v}_0 + \int_0^\infty {\hat
  U} (k, p) e^{-p/t} dp$$ satisfies NS in Fourier space.  Since $\| {\hat
  U} (\cdot, p) \|_{\mu, \beta} < \infty$, it follows that $\| {\hat v}
(\cdot, t) \|_{\mu, \beta} < \infty$ if $\Re ~\frac{1}{t} > \alpha$.  
\end{proof}

\begin{Proposition}[Bounds on $\| {\hat v} (. , t) \|_{\mu+2, \beta}$]
\label{persist}

For the solution ${\hat v} (k, t)$
given in Lemma \ref{localexistence1} for $ t \in [0, \alpha^{-1}]$, 
 we have 
$$ \sup_{t \le T} \| {\hat v} (\cdot, t) \|_{\mu +2, \beta}  < 
C \left ( \| {\hat v}_0 \|_{\mu+2, \beta}, T \right ) < \infty $$  
\end{Proposition}
\begin{proof}
We note from (\ref{nseq}) that if we define $V = \nabla v$, then
${\hat V} = \mathcal{F} [V] = i k {\hat v}$ satisfies, 
\begin{equation}
\label{17}
{\hat V}_t + |k|^2 {\hat V} = - i k \mathcal{P} 
\left [ {\hat v}_j {\hat *} {\hat V}^{[j]} \}
\right ]  + i k {\hat f}  
~~,~~{\hat V}_0 (k) = \mathcal{F} \left [ \nabla v_0 \right ]
\end{equation}
where 
${\hat V}^{[j]} = i k_j {\hat v}_j$.   
Therefore, 
\begin{equation}
\label{17.1}
{\hat V} (k, t) = e^{-|k|^2 t} {\hat V}_0 (k) - i k   
\int_0^t e^{-|k|^2 (t-\tau)}  
\left \{ \mathcal{P} 
\left [ {\hat v}_j {\hat *} {\hat V}^{[j]} \right ] (k, \tau) \}
- {\hat f} (k) \right \} 
\end{equation}  
Therefore, 
\begin{multline}
\label{17.2}
| {\hat V} (\cdot, t) | 
\le 
\frac{e^{-\beta |k|}}{(1+|k|)^\mu} \Bigg \{ 
\| {\hat V}_0 \|_{\mu, \beta} 
+ \\ |k| \int_0^t e^{-|k|^2 (t-\tau) } \left ( \| {\hat f} \|_{\mu, \beta} 
+ 2 C_0 \| {\hat v} (\cdot, \tau) \|_{\mu, \beta}  \|  
\|  {\hat V} (\cdot, \tau) \|_{\mu, \beta} \right ) d\tau \Bigg \}
\end{multline}
Let $ \mathcal{V}_{T_1}$ be the Banach space of continuous
functions $g$ of $k\in \mathbb{R}^3$ and $t \in [0, T_1]$ for which
$$ \| g \|_{T_1} = \sup_{t \in [0, T_1]} \| g (\cdot, t) \|_{\mu, \beta} 
< \infty $$ 
Then, the estimates in (\ref{17.2}), together with the fact
that for any $t \in [0, T]$, $2 C_0 \| {\hat v} (\cdot, t) \|_{\mu, \beta} \le 
{\tilde C} \left ( T, \| {\hat v}_0 \|_{\mu+2,\beta} \right ) $ imply
there exists  $C_1 (T, \|v_0\|_{\mu+2, \beta} )>0$ so that 
\begin{equation}
\label{17.3}
\| {\hat V} \|_{T_1} \le C_1 
\left \{ \sqrt{T_1} \| {\hat V}_1 \|_{T_1} + \| {\hat V}_0 \|_{\mu, \beta}
+ \sqrt{T_1} \|{\hat f} \|_{\mu, \beta}  \right \},
\end{equation}
where we have used the fact that
$$|k| \int_0^{t} e^{-|k|^2 (t-\tau)} d\tau = 
\frac{1 - e^{-|k|^2 t}}{|k|} \le \sqrt{T_1} \sup_{\gamma \in \mathbb{R}^+}
\frac{1 - e^{-\gamma}}{\gamma^{1/2}} \le C_* \sqrt{T_1}, $$ for some $C_*>0$.
Thus, thinking of ${\hat v}$ as given in (\ref{17}), the estimates in
(\ref{17.3}) and similar estimates on ${\hat V}^{[1]} - {\hat V}^{[2]}$ 
show
that for $C_1 \sqrt{T_1} < 1$ the right hand side of (\ref{17}) is contractive
in  $\mathcal{V}_{T_1}$.  We choose $T_1 \le T$.  Therefore, $
\sup_{t \in [0, T_1]} \| {\hat V} (\cdot, t) \|_{\mu, \beta}<\infty$.
Since the choice of $T_1$ depends on $C_1$, which is independent of $\| {\hat
  V}_0 \|_{\mu,\beta}$, we can repeat the same argument in another interval
$[T_1, 2 T_1]$ and so on until we span the whole interval $[0, T]$ over which
$\| {\hat v}_1 (\cdot, t) \|_{\mu, \beta}$ is uniformly bounded.

We can take additional derivative and repeat the same
type argument for $\mathcal{F} [D^2 {\hat v}] = - k k {\hat v}$ to show
that in $\| |k|^2 {\hat v} (\cdot, t) \|_{\mu, \beta}$ is also bounded
uniformly for $t \in [0, T]$. 
In this part of the argument,
we use the prior knowledge that
both $ \| {\hat v} (\cdot, t) \|_{\mu, \beta}$ and
$ \| k {\hat v} (\cdot, t) \|_{\mu, \beta}$ are uniformly bounded in
$[0, T]$ and that 
$$|k|^2 \int_0^{t} e^{-|k|^2 (t-\tau)} \| {\hat f} \|_{\mu, \beta}
d\tau =
\| {\hat f} \|_{\mu, \beta} 
\left (1 - e^{-|k|^2 t} \right ) \le \| {\hat f} \|_{\mu, \beta}
\sup_{\gamma \in \mathbb{R}^+}  
[ 1 - e^{-\gamma} ] \le C \| {\hat f} \|_{\mu, \beta} $$ 
Combining all the results, 
it follows that
$ \| {\hat v} ( \cdot, t) \|_{\mu+2, \beta}$ is bounded for $t \in [0, T]$    
\end{proof}

\noindent{\bf Proof of Theorem \ref{localexistence}.}
This follows from Lemma \ref{localexistence1} and Proposition
\ref{persist}, noting that $ \| {\hat v} (k, t ) \|_{\mu+2, \beta} <
\infty$ implies $ v (x, t) = \mathcal{F}^{-1} [ {\hat v} (\cdot, t) ] (x)\in
{C}^2 (\mathbb{R}^3)$ and so $v$ is a classical solution to (\ref{nseq0}) for
$ \Re ~\frac{1}{t} > \alpha $, which is known to be unique.  From the
definition of $\| \cdot \|_{\mu, \beta}$ it follows that $\| {\hat v}_0
\|_{\mu+2, \beta} < \infty$ and $ \| {\hat f} \|_{\mu, \beta} < \infty$ for
$\beta > 0$ imply $\|{\hat v} (., t)\|_{\mu+2, \beta}<\infty $. 
Thus $v$ preserves the
analyticity strip width for $t \in [0, \frac{1}{\alpha} )$.

\section{Analyticity of ${\hat U} (k, p)$ at $p=0$}

We now consider the case $\beta > 0$.We note that by Remark~\ref{1.6} we can choose
$\mu>3$. 
The starting point of this section
is (\ref{100.3}), which is satisfied by
${\hat W} (k, p) ={\hat U} (k, p) - {\hat v}_1 (k)$. From
Lemma \ref{inteqn}, this is the only solution to (\ref{100.3})
satisfying ${\hat W} (k, 0) =0$.
We seek an potentially alternate  solution to (\ref{100.3}) as a power series, 
\begin{equation}
\label{12.1.0}
{\hat W} (k, p) = \sum_{l=1}^\infty {\hat W}^{[l]} (k) p^l
\end{equation}
Substituting (\ref{12.1.0}) into (\ref{100.3}) and identifying the
coefficients of $p^l,l=0,1$ we get 
\begin{equation}
\label{12.3}
2 {\hat W}^{[1]} = -|k|^2 {\hat v}_1 - i k_j  
P_k
\left [ {\hat v}_{0,j}{\hat *}{\hat v}_1 + 
{\hat v}_{1,j}{\hat *}{\hat v}_0 \right ] , 
\end{equation}
\begin{equation}
\label{12.3.2}
6 {\hat W}^{[2]} = 
-k^2 {\hat W}^{[1]} 
- i k_j P_k 
\left [ 
{\hat v}_{0,j} {\hat *} {\hat W}^{[1]} 
+ {\hat W}^{[1]}_j {\hat *} {\hat v}_0 + 
{\hat v}_{1,j} {\hat *} {\hat v}_1 \right ]
\end{equation}
It follows from (\ref{12.3}) and
Lemma (\ref{lem0.2}) that
\begin{equation}
\label{12.3.1}
\lvert {\hat W}^{[1]} (k, p) \rvert  
\le \frac{e^{-\beta |k|}}{2 (1+|k|)^\mu}
\left ( |k|^2 \| v_1 \|_{\mu, \beta} + 4 C_0 |k| \| v_0 \|_{\mu, \beta}
\| v_1 \|_{\mu, \beta} \right ) 
\end{equation}   
The coefficient of $p^l$ for $l \ge 2$ in (\ref{100.3}) can
be computed as well, using $p^{l_1} * p^{l_2} =p^{l_1+l_2+1} {l_1! l_2!}/{
(l_1+l_2+1)!} $.
Interpreting ${\hat W}^{[0]} = 0$,  we get
\begin{multline}
\label{12.2.5}
(l+1) (l+2) {\hat W}^{[l+1]} = 
-k^2 {\hat W}^{[l]} 
- i k_j P_k  
\left [  \sum_{l_1=1}^{l-2} \frac{l_1! (l-1-l_1)!}{l!} 
{\hat W}_j^{[l_1]}{\hat *} {\hat W}^{[l-1-l_1]}
\right ] \\
- i k_j  P_k 
\left [ 
{\hat v}_{0,j}{\hat *}{\hat W}^{[l]} 
+ {\hat W}^{[l]}_j {\hat *} {\hat v}_0 
+\frac{1}{l} {\hat v}_{1,j} {\hat *} {\hat W}^{[l-1]} 
+ \frac{1}{l} {\hat W}^{[l-1]}_j {\hat *} {\hat v}_1 \right ] 
\end{multline}

\begin{Definition}  
\label{defQ}
It is convenient to define the $n$-th order polynomial
$Q_n$:
$$ Q_n (y) = \sum_{j=0}^n 2^{n-j} \frac{y^j}{j!} $$ 
\end{Definition}

\begin{Lemma} 
\label{lem4}
If $ \| v_0 \|_{\mu+2, \beta} < \infty$, for $\mu > 3$, $\beta > 0$, then
there exist positive constants $A_0,B_0 > 0$ independent of $l$ and $k$ 
so that for any $l \ge 1$ we have
\begin{equation}
\label{12.4}
|{\hat W}^{[l]} (k) | \le e^{-\beta |k|} A_0 B_0^l (1+|k|)^{-\mu} 
\frac{Q_{2l} (|\beta k|)}{(2l+1)^2}
\end{equation}
and
$$ |W^{[l]} (x)| \le \frac{8 \pi A_0 (4 B_0)^l}{(2l+1)^2} ~~,~~
|D W^{[l]} (x)| \le \frac{8 \pi A_0 (4 B_0)^l }{
\beta (2l+1)^2 } ~~,~~
|D^2 W^{[l]} (x)| \le \frac{16 \pi A_0 (4 B_0)^l }{
\beta^2 (2l+1)^2 }
$$
Furthermore, the solution in Lemma \ref{inteqn}, \S 2 
has a convergent series
representation in $p$:
${\hat U} (k, p) = {\hat v}_1 (k) + \sum_{l=1}^\infty {\hat W}^{[l]} 
p^l $ for
$|p| < {(4 B_0)^{-1}}$.
\end{Lemma}

\begin{Remark} Lemma \ref{lem4} is proved by induction on $l$.
For $l=1$, by (\ref{12.3.1})
we just choose
\begin{equation}
\label{12.4.0}
A_0 B_0 \ge 
\frac{18}{\beta^2} \| {\hat v}_1 \|_{\mu, \beta} 
(1 + \beta C_0 \| v_0 \|_{\mu, \beta} )
\end{equation}
Let now $l \ge 2$.
For the induction step, we will estimate 
each term on the right of (\ref{12.2.5}).
\end{Remark}

\begin{Lemma}
\label{lem5.0.0.0}
If for $l \ge 1$, $W^{[l]}$ satisfies (\ref{12.4}), then
$$  \frac{|k|^2 |{\hat W}^{[l]}|}{(l+1)(l+2)} 
\le 
\frac{6 A_0 B_0^{l} e^{-\beta |k|}}{\beta^2 (1+|k|)^\mu}  
\frac{Q_{2l+2} (\beta |k|)}{(2l+3)^2} 
$$  
\end{Lemma}
\begin{proof}
The proof simply follows from the (\ref{12.4}) and
noting that for $y \ge 0$ 
$$ \frac{y^2}{(2l+2)(2l+1)} Q_{2l} (y) 
\le Q_{2l+2} (y)~~~~,~~~~~~
\frac{(2l+3)^2}{(2l+1)(l+2)} \le 3 $$
\end{proof}

\begin{Lemma}
\label{lem5.0}
If $W^{[l]}$ satisfies  (\ref{12.4}), then
for $ l \ge 1$,
\begin{multline*}
\frac{1}{(l+1)(l+2)} \lvert k_j 
P_k 
{\hat u}_{0,j}{\hat *}{\hat W}^{[l]} \rvert 
\le  
2^{\mu} \| v_0 \|_{\mu, \beta}
\frac{9 \pi A_0 B_0^l e^{-\beta |k|}}{\beta^3 
(2l+3)^2 (1+|k|)^\mu} Q_{2l+2} (\beta |k|) \\
\frac{1}{(l+1)(l+2)} \lvert k_j 
P_k
{\hat W}_{j}^{[l]}{\hat *}{\hat u}_{0,j} \rvert 
\le  
2^{\mu} \| v_0 \|_{\mu, \beta}
\frac{9 \pi A_0 B_0^l e^{-\beta |k|}}{\beta^3 
(2l+3)^2 (1+|k|)^\mu} Q_{2l+2} (\beta |k|) 
\end{multline*}
\end{Lemma}
\begin{proof}
We use the estimate (\ref{12.4}) on ${\hat W}^{[l]}$. From
Lemma \ref{lema1.3} for $n=0$, we obtain
\begin{multline*}
\lvert k_j {\hat W}_j^{[l]} {\hat *} {\hat u}_0 \rvert
\le \| v_0 \|_{\mu, \beta} \frac{A_0 B_0^{l}}{(2l+1)^2}
\left ( |k| \int_{k' \in \mathbb{R}^3} 
\frac{e^{-\beta (|k'| + |k-k'|)}}{(1+|k'|)^{\mu} [1+|k-k')]^{\mu}} 
Q_{2l} (\beta |k'|) 
dk' \right ) \\
\le
\frac{\| v_0 \|_{\mu, \beta} A_0 B_0^l}{
(2l+1)^2} 
\sum_{m=0}^{2l} \frac{2^{2l-m}}{m!}   
|k| \int_{k' \in \mathbb{R}^3} 
e^{-\beta (|k'| + |k-k'|)} (1+|k'|)^{-\mu} (1+|k-k'))^{-\mu} |k'|^{2m} 
dk' 
\\
\le \frac{ 2\pi \|v_0 \|_{\mu, \beta} A_0 B_0^l 
2^{\mu} e^{-\beta |k|}}{(2 l+1)^2 
\beta^3 (1+|k|)^\mu} 
\sum_{m=0}^{2l} 2^{2l-m} (m+2) Q_{m+2} (\beta |k|)      
\\
\le \frac{2^{\mu+1} \pi
}{(2 l+1) 
\beta^3 (1+|k|)^\mu} 
 \| v_0 \|_{\mu, \beta} A_0 B_0^l e^{-\beta |k|}
(l+2) Q_{2l+2} (\beta |k|)  
\end{multline*}
The first part of the lemma follows by using (\ref{8.1}) and
checking that
$\frac{2 (2l+3)^2}{(2l+1) (l+1)} \le 9$ for
$l \ge 1$.
The proof of the second part is essentially the same since
$|{\hat W}_j^{[l]}| \le |{\hat W}^{[l]}|$.
\end{proof}

\begin{Lemma}
\label{lem5.0.0}
If $W^{[l-1]}$ satisfies (\ref{12.4}) for any $l \ge 2$, then
\begin{multline*}
\frac{1}{l (l+1)(l+2)} \lvert k_j P_k \left [ 
{\hat u}_{1,j}{\hat *}{\hat W}^{[l-1]} \right ] \rvert 
\le  
2^{\mu} \|v_1 \|_{\mu, \beta} 9 \pi A_0 B_0^{l} (1+|k|)^{-\mu} 
e^{-\beta |k|} \frac{Q_{2l} (\beta |k|) }{\beta^3 (l+2) (2 l +1)^2} \\
\frac{1}{l (l+1)(l+2)} \lvert k_j  P_k
\left [ {\hat W}_{j}^{[l-1]}{\hat *}{\hat u}_{1,j} \right ] \rvert 
\le  
2^{\mu} \| v_1 \|_{\mu, \beta} 9 \pi A_0 B_0^{l-1} (1+|k|)^{-\mu} 
e^{-\beta |k|} \frac{Q_{2l} (\beta |k|) }{\beta^3 (l+2) (2 l +1)^2} 
\end{multline*}
\end{Lemma}
\begin{proof}
The proof is identical to that of Lemma~\ref{lem5.0} with
$l$  replaced by $l-1$ and $v_0$ by $v_1$. 
\end{proof}

\begin{Lemma}
\label{lem5}
If for $l \ge 3$,
$ {\hat W}^{[l_1]}$ and $ {\hat W}^{[l-1-l_1]} $
for $l_1 =1, ..., (l-2)$
satisfy 
(\ref{12.4}), then 
\begin{multline*}
\left\lvert  \frac{k_j}{(l+1)(l+2)} P_k  
\left [  \sum_{l_1=1}^{l-2} \frac{l_1! (l-1-l_1)!}{l!} 
{\hat W}_j^{[l_1]} {\hat *} 
{\hat W}^{[l-1-l_1]}
\right ] \right\rvert \\
\le 
2^\mu 36 A_0^2 B_0^{l-1} (1+|k|)^{-\mu} 
e^{-\beta |k|} \frac{Q_{2l} (\beta |k|)}{\beta^3 (2l+3)^2}
\end{multline*}
\end{Lemma}

\begin{proof}
First note that if we define $l_2=l-1-l_1$, then for $l \ge 3$, 
Lemma \ref{lema1.4} implies
\begin{multline*}
\frac{l_1! l_2!}{l!} 
\left\lvert k_j {\hat W}_j^{[l_1]} \hat{*} {\hat W}^{[l_2]} \right\rvert
\le A_0^2 B_0^{l-1} \frac{(l_1)!(l_2)!}{l! (2l_1+1)^2 (2 l_2+1)^2 } \times \\
|k| \int_{k' \in \mathbb{R}^3} 
e^{-\beta (|k'| + |k-k'|)} (1+|k'|)^{-\mu} (1+|k-k'))^{-\mu} 
Q_{2l_1} (\beta |k'|) Q_{2l_2} (\beta |k-k'|) 
dk' \\
\le \frac{2^{\mu+1} \pi A_0^2 B_0^{l-1} e^{-\beta |k|}}{3 \beta^3 (1+|k|)^\mu} 
\frac{(2 l-1) (2l) (2l+1) l_1 ! l_2!}{l! (2l_1+1)^2 (2 l_2+2)^2 } 
Q_{2l} (\beta |k|) 
\end{multline*}
Therefore, 
\begin{multline*}
\sum_{l_1=1}^{l-2} \frac{l_1! l_2!}{l! (l+1)(l+2)} 
\lvert k_j {\hat W}_j^{[l_1]} \hat{*} {\hat W}^{[l_2]} \rvert \\
\le 
\frac{2^{\mu+2} \pi e^{-\beta |k|} Q_{2l} (\beta |k|) (2l-1) (2l+1)}{ 
3 (l+1) (l+2) \beta^3 (1+|k|)^\mu} 
\sum_{l_1=1}^{l-2}
\frac{l_1 ! l_2 !}{(l-1)! (2l_1+1)^2 (2l_2+1)^2 } 
\end{multline*}
and the proof follows noting that $\frac{l_1! l_2!}{(l-1)!} = \frac{l_1!
  l_2!}{(l_1+l_2)!} \le 1$ and checking $\frac{4 (2l-1) (2l+1)}{(l+1) (l+2)}
\le 16$; by breaking up the sum in the ranges: $ l_1 \le (l-1)/2$ and $l_1 >
(l-1)/2 $ (in which $l_2 \le (l-1)/2$) it is easily seen that for some
$C_*>0$ and any $l \ge 3$ we have
$$ \sum_{l_1=1}^{l-2}
\frac{1}{(2l_1+1)^2 (2l_2+1)^2} \le \frac{C_*}{(2l+3)^2}, $$
where $C_* = 1.07555\cdots $ (the upper-bound being
achieved at $l=4$). 
\end{proof}

\begin{Lemma}
\begin{equation}
| {\hat W}^{[2]} | \le  \frac{e^{-\beta |k|}}{(1+|k|)^\mu} 
\frac{Q_{4} (\beta |k|)}{7^2} \left ( 
\frac{A_0 B_0}{\beta^2}
+ A_0 B_0 \| v_0 \|_{\mu, \beta} \frac{2^\mu 36 \pi}{\beta^2} 
+ \| v_1 \|_{\mu, \beta}^2 \right )   
\end{equation}
and therefore 
${\hat W}^{[2]}$ satisfies  (\ref{12.4}) if
\begin{equation}
\label{ab2cond}
A_0 B_0^2 \ge  
\frac{3 A_0 B_0}{\beta^2}
+ A_0 B_0 \| v_0 \|_{\mu, \beta} \frac{2^\mu 36 \pi}{\beta^2} 
+ \frac{C_0}{\beta} \| v_1 \|_{\mu, \beta}^2 
\end{equation}
\end{Lemma} 
\begin{proof}
We use Lemmas \ref{lem5.0.0.0}, \ref{lem5.0} and 
\ref{lem0.1} to estimate different terms on the right hand side
of (\ref{12.2.5}) for $l=1$.  
\end{proof}

\noindent{\bf Proof of Lemma \ref{lem4}}

We use Lemmas \ref{lem5.0.0.0}, \ref{lem5.0}
\ref{lem5.0.0} and \ref{lem5} to estimates the
terms on the right hand side of (\ref{12.2.5}) and note that
$Q_{2l} (y) \le \frac{1}{4} Q_{2l+2} (y)$. Hence, combining
all the estimates,
we obtain for $l \ge 2$,
\begin{multline*}
| {\hat W}^{[l+1]} | \le 
A_0 B_0^{l-1} 
\frac{Q_{2l+2} (\beta|k|)e^{-\beta |k|}}{(2l+3)^2 (1+|k|)^\mu}\\ \times 
\left \{ \frac{6}{\beta^2} B_0 + 2^{\mu} \frac{18 \pi}{\beta^3} 
B_0 \| v_0 \|_{\mu, \beta} + 
\frac{2^{\mu} 18\pi (2l+3)^2 }{(l+2) (2l+1)^2 \beta^3} 
\|v_1 \|_{\mu, \beta} + 
\frac{9 A_0 2^\mu}{\beta^3} \right \}     
 \\ \le  
\frac{A_0 B_0^{l+1} e^{-\beta |k|}}{(1+|k|)^\mu (2l+3)^2} Q_{2l+2} (\beta|k|)
\end{multline*}
for large enough $B_0$ so that
\begin{equation}
\label{Bcond}
\left \{ \frac{6}{\beta^2} B_0 + 2^{\mu} \frac{18\pi}{\beta^3} 
B_0 \| v_0 \|_{\mu, \beta} + 
\frac{2^{\mu} 18 \pi }{\beta^3} \|v_1 \|_{\mu, \beta} + 
\frac{9 A_0 2^\mu}{\beta^3} \right \} \le B_0^2 
\end{equation}
Combining (\ref{Bcond}) with (\ref{12.4.0}) and (\ref{ab2cond}), we that
(\ref{12.4}) is satisfied for any $l \ge 1$.  Therefore, it follows that
$\sum_{l=1}^\infty {\hat W}^{[l]} (k) p^l $ is convergent for $ |p| <
\frac{1}{4 B_0}$. The recurrence relations (\ref{12.3}),(\ref{12.3.2}) and
(\ref{12.2.5}) imply that $\sum_{l=1}^\infty {\hat W}^{[l]} (k) p^l $ is
indeed a solution to (\ref{100.3}), which is zero at $p=0$. However, 
from \S2 
Lemma \ref{inteqn},
we know that there is a unique ${\hat W} = {\hat U} (k, p) - {\hat v}_1 (k)$
with this property in $\mathcal{A}_L^\infty$, which for
sufficiently small $L$ includes analytic functions at the origin. Therefore
$$ {\hat U} (k, p) = {\hat v}_1 (k) + \sum_{l=1}^\infty {\hat W}^{[l]} (k)
p^l $$
Moreover, from the well-known relation between a function and
its Fourier transform,
$ \| f \|_{L^\infty (\mathbb{R}^3)} \le 
\| {\hat f} \|_{L^1(\mathbb{R}^3}$,
the inequalities involving $W^{[l]} (x)$ and its $x$-derivatives follow.

\section{Estimates on $\partial_p^l {\hat W} (k, p)$ and
proof of Theorem \ref{T2}}

In this section, we find inductively (in $l$) that $ {\hat
  W}^{[l]}:=\partial_p^l \hat {W}/l!$ exists for any $l$ and ${\hat W}^{[l]}$
generate power series (\ref{100}) with $p_0-$ independent radius of
convergence. This does not necessarily imply in itself that the series
converges to $\hat {W}$.  The fact that these objects do coincide locally will
be shown in Lemma~\ref{lemancont}. This leads to proof of Theorem \ref{T2}.

\begin{Definition}
It is convenient to define for $l \ge 1$, 
$$ {\hat W}^{[l]} (k, p) = \frac{1}{l!} \partial^l_p {\hat W} (x, p) ~~,~~
$$
It is also convenient to define ${\hat W}^{[0]} (k, p) = {\hat W} (k, p)
= {\hat U} (k, p) - {\hat v}_1 (k) $.
\end{Definition}

The proof  therefore reduces to
finding appropriate bounds on ${\hat W}^{[l]} (k, p)$. 
The main result proved in this section is the Lemma \ref{lem40},
which, using Lemma \ref{lemancont},
leads directly to the proof of Theorem \ref{T2}.

Proposition \ref{inteqnA} implies that ${\hat U} \in 
\mathcal{A}^{\alpha'}$ for $\alpha' > \alpha$, with
$\alpha$ chosen large enough to satisfy (\ref{ensure2}). In
particular, if we choose $\alpha' = \alpha+1$, it follows that
${\hat W}^{[0]} (k, p) = {\hat U} (k, p) - {\hat v}_k (k) $
satisfies
\begin{equation}
\label{W0inequal}
|{\hat W}^{[0]} (k, p)| \le 
\frac{3 e^{-\beta |k|+\alpha' p} \| {\hat v}_1 \|_{\mu, \beta} }{
(1+p^2) (1+|k|)^\mu } 
\end{equation}
In the rest of this section, with some
abuse of notation, we will replace $\alpha'$ by $\alpha$.

\begin{Lemma} 
\label{lem40}
If $\| v_0 \|_{\mu+2, \beta} < \infty$, $\mu > 3$,
there exists positive constants $A$, $B$ independent of $l$, $k$ and $p$ 
so that
for any $l \ge 0$
\begin{equation}
\label{104}
|{\hat W}^{[l]} (k, p) | \le \frac{
e^{\alpha p} e^{-\beta |k|}}{(1+p^2) (1+|k|)^{\mu}} A B^l  
\frac{Q_{2l} (|\beta k|)}{(2l+1)^2}
\end{equation}
The series (\ref{100})
converges uniformly for any $p_0 \ge 0$ for $|p-p_0| < \frac{1}{4B}$.
\end{Lemma}

\begin{Remark} The proof requires some further
lemmas. We will use induction on $l$. 
Clearly, from (\ref{W0inequal}), the conclusion
is valid for $l=0$, when
\begin{equation}
\label{conditionA}
A = 3 \| {\hat v}_1 \|_{\mu, \beta} 
\end{equation} 
We assume (\ref{104}) for 
$l \ge 0  $ and then establish it for $l+1$.
We obtain a recurrence relation for 
${\hat W}^{[l+1]} (k, .)$ for any $k \in \mathbb{R}^3$  
in terms of ${\hat W}^{[j]} (k , .)$ for $j \le l$. 
\end{Remark}

Taking $\partial_p^l $ in
(\ref{100.3}) and dividing by $l!$, we obtain
\begin{multline}
\label{102}
p \partial_p^2 {\hat W}^{[l]} + (l+2) \partial_p {\hat W}^{[l]} + 
|k|^2 {\hat W}^{[l]} = \\
-i k_j P_k 
\left [ \int_0^{p} \left \{ {\hat W}_j^{[l]} (\cdot, p-s) 
{\hat *} {\hat W}^{[0]} (\cdot, s) \right \} ds + \sum_{l_1=1}^{l-1} 
\frac{l_1! (l-1-l_1)!}{l!} {\hat W}_j^{[l_1]} (\cdot, 0) {\hat *}
{\hat W}^{[l-1-l_1]} (\cdot, p) \right ]  
\\
-i k_j P_k  
\left [ 
{\hat v}_{0,j}\hat{*}{\hat W}^{[l]} +    
{\hat v}_{0}\hat{*}{\hat W}_j^{[l]} +    
\frac{1}{l} {\hat v}_{1,j}{\hat *}{\hat W}^{[l-1]} 
+ \frac{1}{l} {\hat v}_{1}{\hat *}{\hat W}_j^{[l-1]} + 
{\hat v}_{1,j}{\hat *}{\hat v}_1
\delta_{l,1} \right ] \\
-i k_j P_k \left [ {\hat v}_{0,j} {\hat *} {\hat v}_1 
+ {\hat v}_{1,j} {\hat *} {\hat v}_0 \right ] \delta_{l,0} 
\equiv {\hat R}^{(l)} (k, p) 
\end{multline}

\begin{Lemma} 
\label{URbound}
For any $l \ge 0$, for some absolute constant $C_6>0$, if ${\hat W}^{[l]} (k,
p)$ satisfies (\ref{102}) and is bounded at $p=0$, then 
${\hat W}^{[l+1]} (k, p)$ is bounded in terms of 
${\hat
R}^{(l)} (k, p)$, defined in (\ref{102}):
$$ |{\hat W}^{[l+1]} (k, p) \le \frac{C_6}{(l+1)^{5/3}} 
\sup_{p' \in [0, p]} |{\hat R}^{[l]} (k, p)| 
+ \frac{|k|^2 |{\hat W}^{[l]} (k, 0)|}{(l+1)(l+2)}  $$ 
\end{Lemma}

\begin{proof}
We invert 
the operator on the left hand side of (\ref{102}). With 
the requirement
that ${\hat W}^{[l]} $ is bounded at $p=0$, we obtain
\begin{multline}
\label{202}
{\hat W}^{[l]} (k, p) = \int_0^p 
\mathcal{Q} \left (z (p), 2 |k| \sqrt{p'} \right ) {\hat R}^{(l)} 
(k, p') dp' \\+ 2^{(l+1)} (l+1)! {\hat W}^{[l]} (k, 0) 
\frac{J_{l+1} (z)}{z^{l+1}} ~~,~{\rm where}
~z = 2 |k| \sqrt{p}  
\end{multline}
and
\begin{multline}
\label{203}
\mathcal{Q} (z, z') =  
\pi z^{-(l+1)} \left [ -J_{l+1} (z) {z'}^{(l+1)}    
Y_{l+1} (z') + {z'}^{(l+1)} J_{l+1} (z') Y_{l+1} (z) 
\right ]  
\end{multline}
On taking the first derivative with respect to $p$, we 
obtain
\begin{multline}
\label{204}
(l+1) {\hat W}^{[l+1]} (k, p)\\ = 
\frac{|k|}{\sqrt{p}} \int_0^p 
\mathcal{Q}_z \left (2 |k| \sqrt{p}, 2 |k| \sqrt{p'} \right ) 
{\hat R}^{(l)} (k, p') dp' - 2^{l+2} (l+1)! |k|^2 
\frac{J_{l+2} (z)}{z^{l+2}}
{\hat W}^{[l]} (k, 0)   
\end{multline}
Using again the properties of Bessel functions \cite{Abramowitz} we get
\begin{multline}
\label{204.1}
\frac{1}{z} \mathcal{Q}_z (z, z')
= 
\frac{\pi}{z} 
\left [ -\left (\frac{J_{l+1} (z)}{z^{l+1}} \right )^\prime  
{z'}^{(l+1)}    
Y_{l+1} (z') + {z'}^{(l+1)} J_{l+1} (z') 
\left ( \frac{Y_{l+1} (z)}{z^{l+1}} \right )^\prime \right ] \\
= 
\pi
\left [ \frac{J_{l+2} (z)}{z^{l+2}}   
{z'}^{(l+1)}    
Y_{l+1} (z') - {z'}^{(l+1)} J_{l+1} (z') 
\frac{Y_{l+2} (z)}{z^{l+2}} \right ] 
\end{multline}
It is also known \cite{Abramowitz} that
$$ 2^{l+2} (l+1)!  
\Big | \frac{J_{l+2} (z)}{z^{l+2}} \Big | \le \frac{1}{(l+2)} 
$$
Using (\ref{204.1}) and the known uniform
asymptotics of Bessel functions for
large $l$
\cite{Abramowitz}, it is easily to see
that $C_*$
independent of  $l$ so that
$$ \int_0^z \frac{z'}{z} | Q_z (z, z') | dz'  \le \frac{C_*}{(l+1)^{2/3}} $$  
It follows that 
\begin{equation}
\label{205}
(l+1) | {\hat W}^{[l+1]} (k, p) | \le \sup_{p' \in [0,p]}   
| {\hat R}^{(l)} (k, p') | \int_0^z \frac{z'}{z} | Q_z (z, z') | dz'  
+ \frac{|k|^2}{(l+2)} | {\hat W}^{[l]} (k, 0) |
\end{equation}
Therefore, it follows that 
\begin{equation}
\label{hatURbound}
\Big | {\hat W}^{[l+1]} (k, p) \Big | \le \frac{C_6}{(l+1)^{5/3}}  
\sup_{p' \in [0,p]}   
| {\hat R}^{(l)} (k, p') | + 
\frac{|k|^2 |{\hat W}^{[l]} (k, 0)|}{(l+1)(l+2)}   
\end{equation}
\end{proof}

\begin{Remark}
We now find bounds on the different terms in
${\hat R}^{(l)} (k, p)$.
\end{Remark}

\begin{Lemma}
\label{lem50.0}
If $W^{[l]}$ satisfies (\ref{104}), for $l \ge 0$ then
\begin{multline*}
\lvert k_j 
P_k 
\left ( {\hat v}_{0,j}{\hat *}{\hat W}^{[l]} \right ) \rvert 
\le  
C_1 \| {\hat v}_0 \|_{\mu, \beta}
\frac{(l+1)^{2/3} A B^l e^{-\beta |k|+\alpha p} Q_{2l+2} (\beta |k|) }{ 
(2l+1) (1+|k|)^\mu (1+p^2)} \\
\lvert k_j P_k \left (  
{\hat W}_{j}^{[l]}{\hat *}{\hat v}_{0,j} \right ) \rvert 
\le  
C_1 \| {\hat v}_0 \|_{\mu, \beta}
\frac{(l+1)^{2/3} A B^l e^{-\beta |k|+\alpha p} Q_{2l+2} (\beta |k|) }{
(2l+1) (1+|k|)^\mu (1+p^2)} 
\end{multline*}
\end{Lemma}
\begin{proof}
We use (\ref{104}). From
Lemma \ref{lemaQ2l}, we obtain 
\begin{multline*}
(1+p^2) e^{-\alpha p} \lvert k_j {\hat W}_j^{[l]}{\hat *}{\hat v}_0 \rvert
\le \| {\hat v}_0 \|_{\mu, \beta} \frac{A B^{l} }{(2l+1)}
|k| \int_{k' \in \mathbb{R}^3} 
\frac{e^{-\beta (|k'| + |k-k'|)}}{(1+|k'|)^{\mu} [1+|k-k')]^{\mu}} 
Q_{2l} (\beta |k'|) 
dk' \\
\le C_1 (l+1)^{2/3} \| {\hat v}_0 \|_{\mu, \beta} \frac{A B^{l} }{(2l+1)}
\frac{e^{-\beta |k|}}{(1+|k|)^\mu} Q_{2l+2} (\beta |k|)
\end{multline*}
The first part of the Lemma follows.
The proof of the second part is essentially the same since
$|{\hat W}_j^{[l]}| \le |{\hat W}^{[l]}|$.
\end{proof}

\begin{Lemma}
\label{lem50.0.0}
If $W^{[l-1]}$ satisfies (\ref{104}) for $l \ge 1$, then
\begin{multline*}
\Big \lvert \frac{k_j}{l} 
P_k \left [
{\hat v}_{1,j}{\hat *} {\hat W}^{[l-1]} \right ] \Big \rvert 
\le  
C_1 \|v_1 \|_{\mu, \beta} A B^{l-1} 
\frac{e^{-\beta |k| + \alpha p}}{(1+p^2)(1+|k|)^\mu} 
\frac{l^{2/3} Q_{2l} (\beta |k|) }{l (2 l -1) } \\
\Big \lvert \frac{k_j}{l} 
\left (1-\frac{k (k \cdot )}{|k|^2} \right ) 
{\hat v}_{1}\hat{*}{\hat W}_j^{[l-1]} \Big \rvert 
\le  
C_1 \|v_1 \|_{\mu, \beta} A B^{l-1} 
\frac{e^{-\beta |k| + \alpha p}}{(1+p^2)(1+|k|)^\mu} 
\frac{l^{2/3} Q_{2l} (\beta |k|) }{l (2 l -1)} 
\end{multline*}
\end{Lemma}
\begin{proof}
The proof is identical to Lemma \ref{lem50.0}
replacing
$l$ by $l-1$ and ${\hat v}_0$ by ${\hat v}_1$. 
\end{proof}

\begin{Lemma}
\label{lem50.0.0.new}
If $W^{[l]}$ satisfies (\ref{104}), then for $l \ge 1$,
\begin{equation*}
\Big \lvert \frac{k_j}{l} 
P_k \left [ 
{\hat W}^{[l-1]} (\cdot, 0) {\hat *} {\hat W}^{[0]} (\cdot,p) \right ] 
\Big \rvert 
\le  
C_1 
\frac{(l+1)^{2/3} A^2 B^{l-1} e^{-\beta |k|+\alpha p} Q_{2l} (\beta |k|) }{ 
l (2l-1) (1+|k|)^\mu (1+p^2)} 
\end{equation*}
\end{Lemma}
\begin{proof}
Noting that
$$ | {\hat W}^{[0]} (k, p) | \le A
\frac{e^{-\beta |k| + \alpha p}}{(1+|k|)^\mu (1+p^2)} $$
and   
$$ | {\hat W}^{[l-1]} (k,0) 
\le 
\frac{e^{-\beta |k|} }{(2l-1)^2 (1+|k|)^\mu} A B^{l-1} Q_{2l-2} (\beta |k|) $$
the rest of the proof is very similar to the proof of Lemma
\ref{lem50.0}
\end{proof}

\begin{Lemma}
\label{lem50}
If ${\hat W}^{[l_1]}$ and ${\hat W}^{[l-1-l_1]}$ for $l_1 =1, ..(l-2)$
for $l \ge 2$
satisfy 
(\ref{104}), then 
\begin{multline*}
\Big \lvert  k_j  
P_k
\left [  \sum_{l_1=1}^{l-2} \frac{l_1! (l-1-l_1)!}{l!} 
{\hat W}_j^{[l_1]} (\cdot, 0) {\hat *}
{\hat W}^{[l-1-l_1]} (\cdot, p)
\right ] \Big \rvert \\
\le 
C_8 2^{\mu +1} \pi A^2 B^{l-1} \frac{e^{-\beta |k| + \alpha p}}{3 \beta^3 
(1+p^2) (1+|k|)^{-\mu}} 
\frac{l Q_{2l} (\beta |k|)}{(2l+3)^2} ~~;~~{\rm where}~C_8 = 82 
\end{multline*}
\end{Lemma}

\begin{proof}
First note that if we define $l_2=l-1-l_1$, then for $l \ge 2$, using
Lemma \ref{lema1.4}, we get
\begin{multline*}
\frac{l_1! l_2!}{l!} 
\lvert k_j {\hat W}_j^{[l_1]} (\cdot, 0) {\hat *} {\hat W}^{[l_2]} (\cdot, p) \rvert
\le \frac{e^{\alpha p}}{(1+p^2)} 
A^2 B^{l-1} \frac{l_1 ! l_2 !}{l! (2l_1+1)^2 (2l_2+1)^2}  \\
\times |k| \int_{k' \in \mathbb{R}^3} 
e^{-\beta (|k'| + |k-k'|)} (1+|k'|)^{-\mu} (1+|k-k'))^{-\mu} 
Q_{2l_1} (\beta |k'|) Q_{2l_2} (\beta |k-k'|) 
dk' \\
\le \frac{A^2 B^{l-1} 2^{\mu+1} \pi e^{-\beta |k|+\alpha p}}{
3 \beta^3 (1+p^2) (1+|k|)^\mu} 
\frac{l_1! l_2! (2l) (2l-1) (2l+1)}{l! (2l_1+1)^2 (2l_2+1)^2} 
Q_{2l} (\beta |k|) 
\end{multline*}
Therefore, 
\begin{multline*}
\sum_{l_1=1}^{l-2} \frac{l_1! l_2!}{l!} 
\lvert k_j {\hat W}_j^{[l_1]} (.; 0) \hat{*} {\hat W}^{[l_2]} (.; p)\rvert
\\
\le 
\frac{2^{\mu} A^2 B^{l} e^{-\beta |k|+\alpha p} 
Q_{2l} (\beta |k|)}{\beta^3 (1+|k|)^\mu (1+p^2)} \sum_{l_1=1}^{l-2}
\frac{l_1! l_2! (2 l) (2l+1) (2l-1)}{l! (2l_1+1)^2 (2l_2+1)^2} 
\end{multline*}

We claim that for $l \ge 2$,
with $l_1 \ge 1$, $l_2 = l - l_1 - 1 \ge 1$,  
$$ 
\sum_{l_1=1}^{l-1} \frac{l_1! l_2! (2l ) (2l-1) (2l+1)}{l! (l-1) 
(2l_1+1)^2 (2l_2+1)^2}   
\le \frac{C_8 l}{(2l +3)^2}  
$$
for some 
$C_8 $ independent of $l$; $C_8$ is bounded by 82.

Proving the above bound
only requires consideration for 
sufficiently large $l$. We will therefore assume $l \ge 5$.
Further, consider summation terms other than
$l_1=1$ and $l_2 = 1$. So, we may assume $l_1, l_2 \ge 2$.
Then, we claim that
\begin{multline}
\label{L1L2}
\frac{l_1! l_2! 2l (2l+1) (2l-1)}{
l! (2l_1+1)^2 (2l_2+1)^2}   
= \left (\frac{(l_1-2)! (l_2-2)!}{(l-5)!} \right ) 
\left ( \frac{l_1 (l_1-1) l_2 (l_2-1)}{(2l_1+1)^2 (2l_2+2)^2} \right ) \times
\\
\left ( \frac{2l (2l-1) (2l+1)}{l (l-1) (l-2) (l-3) (l-4)} \right ) 
\le \frac{12}{(2l + 3)^2} 
\end{multline}
This follows since the first two parenthesis term on the right of
(\ref{L1L2}) is clearly bounded, while the last term is
a cubic in $l$ divided by fifth order polynomial, and simple estimates
give the upperbound of 12.
Therefore, for $l \ge 5$,
$$ \sum_{l_1=2}^{l-3}
\frac{l_1! l_2! 2l (2l+1) (2l-1)}{
l! (2l_1+1)^2 (2l_2+1)^2}  \le 12 \frac{(l-4)}{(2l+3)^2}   
$$
For $l_1=1$ or $l_2=1$, clearly 
$$
\frac{l_1! l_2! 2l (2l+1) (2l-1)}{
l! (2l_1+1)^2 (2l_2+1)^2} =    
\frac{ (l-2) ! 2l (2l+1) (2l-1)}{
9 l! (2 l -3)^2} =    
\frac{ 2l (2l+1) (2l-1)}{
9 l (l-1) (2 l -3)^2} \le 82 \frac{l}{(2l+3)^2}    
$$
\end{proof}

\begin{Lemma}
\label{lem50.0.0.0}
If ${\hat W}^{[l]}$ 
satisfies 
(\ref{104}), then for $l \ge 0$, 
\begin{multline*}
\left\lvert  k_j  
\left (1-\frac{k (k \cdot )}{|k|^2} \right )
\int_0^{p} {\hat W}_j^{[l]} (.; p-s){\hat *} 
{\hat W}^{[0]} (.; s) ds
\right\rvert \\
\le 
C_1 (l+1)^{2/3} 
A^2 B^l \frac{e^{-\beta |k| +\alpha p} Q_{2l+2} (\beta |k|)}{ 
(1+p^2) (1+|k|)^{\mu} (2l+1)} 
\end{multline*}
\end{Lemma}
\begin{proof}
We note that Lemma \ref{lemaQ2l} implies
\begin{multline*}
\left\lvert k_j 
\int_{k' \in \mathbb{R}^3} \int_0^{p} {\hat W}_j^{[l]} (k', p-s)  
{\hat W}^{[0]} (k-k'; s) ds dk' \right\rvert 
\le \frac{ A^2 B^l e^{\alpha p}}{(1+p^2) (2l+1)^2} 
\\
\times |k'| \int_{k' \in \mathbb{R}^3} 
\frac{e^{-\beta |k'|-\beta |k-k'|}}{(1+|k'|)^\mu (1+|k-k')^\mu}
Q_{2l} (\beta k')  dk'
\\ 
\le \frac{ C_1 (l+1)^{2/3} A^2 B^l e^{\alpha p-\beta |k|}}{(2l+1) (1+p^2) 
(1+|k|)^\mu} Q_{2l+2} (\beta |k|)  
\end{multline*}
\end{proof}

\begin{Lemma}
\label{loneterm}
$$ \left\lvert k_j \left (1 - \frac{k (k \cdot .)}{|k|^2} \right ) 
\left ( {\hat v}_{0,j} {\hat *} {\hat v}_1  + {\hat v}_{1,j} {\hat *} 
{\hat v}_0 \right )
\right\rvert  
\le 4 C_0 |k| \frac{e^{-\beta |k|}}{(1+|k|)^\mu} \| {\hat v}_0 \|_{\mu, \beta}
\| {\hat v}_1 \|_{\mu, \beta} $$
$$ \left\lvert k_j \left (1 - \frac{k (k \cdot .)}{|k|^2} \right ) 
{\hat v}_{1,j} {\hat *} {\hat v}_1  \right\rvert  
\le 2 C_0 |k| \frac{e^{-\beta |k|}}{(1+|k|)^\mu} 
\| {\hat v}_1 \|^2_{\mu, \beta} $$
\end{Lemma}
\begin{proof}
This follows simply from the observation that  
$$ \left\lvert k_j \left (1 - \frac{k (k \cdot .)}{|k|^2} \right ) 
{\hat v}_{0,j} {\hat *} {\hat v}_1  \right\rvert 
\le 2 |k| \| {\hat v}_1 \|_{\mu, \beta} \|{\hat v}_0 \|_{\mu, \beta}
\int_{k' \in \mathbb{R}^3} 
\frac{e^{-\beta |k'| - \beta |k-k'|}}{(1+|k'|)^\mu (1+|k-k'|)^\mu } dk' 
$$
and using (\ref{eqlem0.1}) to bound the
convolution. Other parts of the Lemma follow similarly.
\end{proof}

\begin{Lemma}
\label{W2bound}
$$ \lvert {\hat W}^{[1]} (\cdot, p) \rvert 
\le \frac{e^{-\beta |{ k}|+\alpha p}}{(1+|k|)^\mu (1+p^2)} A B 
Q_2 (\beta |k|)
$$
with 
\begin{equation}
\label{ABbound}
AB \ge 
\left ( 2 C_1 \| v_0 \|_{\mu, \beta} A + C_1 A^2 + 
\frac{2 C_0}{\beta} \| {\hat v}_0 \|_{\mu, \beta} 
\| v_1 \|_{\mu, \beta} \right ) 
\end{equation}
\end{Lemma}
\begin{proof}
Combining Lemmas \ref{lem50.0}, \ref{lem50.0.0.0} and \ref{loneterm} with
 (\ref{hatURbound}) for $l=0$, we obtain 
\begin{multline*}
\lvert {\hat W}^{[1]} (\cdot, p) \rvert 
\le \frac{e^{-\beta |{k}|+\alpha p}}{(1+|k|)^\mu (1+p^2)} Q_2 (\beta |k|)
\Bigg ( 2 C_1 \| v_0 \|_{\mu, \beta} A + C_1 A^2  + 
\frac{2 C_0}{\beta} \| {\hat v}_0 \|_{\mu, \beta} 
\| {\hat v}_1 \|_{\mu, \beta} \Bigg )  \\ 
\le \frac{e^{-\beta |{ k}|+\alpha p}}{(1+|k|)^\mu (1+p^2)} A B Q_2 (\beta |k|)
\end{multline*}
\end{proof}

\noindent{\bf Proof of Lemma \ref{lem40}}

From Lemmas \ref{lem50.0}-\ref{lem50.0.0.0} and
\ref{loneterm} (the latter is only needed for $l=1$), 
it follows that ${\hat R}^{(l)}$ (cf.
(\ref{104})) satisfies
\begin{multline*}
|{\hat R}^{(l)} | \le A B^{l-1} 
\frac{e^{-\beta |k| + \alpha p}}{(2l+3)^2 (1+p^2) (1+|k|)^\mu}  
Q_{2l+2} (\beta |k|) \\       
\times \Bigg [ A B C_1 \frac{(l+1)^{2/3} (2l+3)^2 }{(2l+1)} +      
\frac{A C_1 (l+1)^{2/3} (2l+3)^3}{ 4 l (2l -1)} + 
\frac{82~2^{\mu+1} \pi A l}{
12 \beta^3 } + \frac{2 C_1 \| {\hat v}_1 \|_{\mu, \beta} (2l+3)^3}{4 l^{1/3}
(2 l -1) } \\
+ \frac{2 C_0 (l+1)^{2/3} (2l+3)^2 \| {\hat v}_0 \|_{\mu, \beta} }{ (2 l+1)}
+ \frac{25 C_0}{\beta} (1+p^2) 
e^{-\alpha p} \|{\hat v}_1 \|_{\mu, \beta}^2 \delta_{l,1}  
\Bigg ]
\end{multline*}
Noting that 
$e^{-\alpha p} (1+p^2) \le 1 $ and
$$\sup_{p' \in [0,p]} \frac{e^{\alpha p'}}{1+{p'}^2} =
\frac{e^{\alpha p}}{1+p^2} $$ 
for $\alpha \ge 1$, 
it follows from 
Lemma \ref{URbound} and the above bounds that (\ref{104})
holds when $l$ is replaced by $l+1$, provided
$B$ is chosen large enough 
to satisfy (\ref{ABbound}) and
\begin{multline}
\label{Bfcond}
C_6 \Bigg [ A B C_1 \frac{(2l+3)^2 }{(l+1) (2l+1)} +      
\frac{A C_1 (2l+3)^3}{ 4 l (l+1) (2l -1)} + 
\frac{(82) 2^{\mu+1} \pi A l}{
12 \beta^3 (l+1)^{5/3}} + 
\frac{2 C_1 \| {\hat v}_1 \|_{\mu, \beta} (2l+3)^3}{4 l^{1/3}
(l+1)^{5/3} (2 l -1) } \\
+ \frac{2 C_0 (2l+3)^2 \| {\hat v}_0 \|_{\mu, \beta} }{
(l+1)(2 l+1)}
+ \frac{25 C_0}{\beta} (1+p^2) 
e^{-\alpha p} \|{\hat v}_1 \|_{\mu, \beta}^2 \delta_{l,1}  
\Bigg ] + \frac{100 B}{9 \beta^2} \le B^2,   
\end{multline}
for any $l \ge 1$,
with $A$ given by (\ref{conditionA}).
From the
asymptotic behavior of the left hand side of (\ref{Bfcond})
as $l \rightarrow \infty$ and recalling that constants $C_0$, $C_1$ and
$C_6$ are independent of $l$, it follows that 
$B$ can be chosen independent of $l$.
Therefore, by induction, (\ref{104}) follows for all $l$. The proof
of Lemma \ref{lem40} is complete.

From (\ref{104}), after noting that
that $Q_{2l} (|q|) \le 4^l e^{-|q|/2}$, it follows that  
\begin{equation}
\label{100}
{\tilde W} (k, p; p_0) =
\sum_{l=0}^\infty {\hat W}^{[l]} (k, p_0) (p-p_0)^l := {\hat W}_1 (k, p)
\end{equation}
is convergent for $|p-p_0| < \frac{1}{4B}$ for $B$ independent of 
$p_0 \in \mathbb{R}^+$. The following Lemma shows that
${\tilde W} (k, p; p_0)$ is indeed the local representation of
the solution
${\hat W} (k, p)$ to (\ref{100.3}).

\begin{Lemma}
\label{lemancont}
The unique solution to (\ref{100.3}) satisfying ${\hat W} (k, 0) =0$,  
given by ${\hat W} (k, p) = {\hat U} (k, p) - {\hat v}_1 (k)$,
where ${\hat U} (k, p)$ is determined in \S 2  in Lemma \ref{inteqn},
has the local representation ${\tilde W} (k, p; p_0)$ in a neighborhood
of $p_0 \in \mathbb{R}^+$. Therefore, ${\hat W} (k, .)$ 
(and therefore ${\hat U} (k, .) $) is analytic in $\mathbb{R}^+ \cup \{0 \}$. 
\end{Lemma}

\begin{proof} 
  First, by permanence of relations (for analyticity of
  convolutions, see e.g., \cite{Duke}), it follows
  that if $\hat{V}$ is an analytic
  solution of an equation of the form (\ref{100.3}) on an interval $[0,L]$ and
  $\hat{V}$ has analytic continuation on $[0,L']$ with $L'>L$, then the
  equation is automatically satisfied in the larger interval. Therefore, if we analytically continue $\hat{W}$ to 
$\mathbb{R}^+$, the analytic continuation will
automatically satisfy (\ref{100.3}) and will therefore be the same
as ${\hat W} (k, p)$.

From \S3, Lemma \ref{lem4}, we know that 
the actual solution to (\ref{100.3}) satisfying ${\hat W} (k, 0)=0$, 
is unique, and given by
$$ {\hat W} (k, p) = {\tilde W} (k, p; 0)$$  
for $ |p| < (4B)^{-1}$.

We now choose a sequence of $\left \{ p_{0,j} \right \}_{j=0}^\infty$,
with $p_{0,j} = j/(8B)$ and define the intervals 
$\mathcal{I}_j = \left ( p_{0, j} - 1/(4B), p_{0, j} +1/(4B)\right )$.
Consider the sequence of analytic functions 
$\left \{ {\tilde W} (k, p; p_{0, j} ) \right \}_{j=0}^\infty$.
Since $p_{0, 1} \in \mathcal{I}_0 \cap \mathcal{I}_1$, it follows from
 (\ref{100}) that ${\hat W} (k, p) $ has analytic continuation 
to $\mathcal{I}_1$, namely ${\tilde W} (k, p; p_{0,1})$.
Again $p_{0,2} \in \mathcal{I}_1 \cap \mathcal{I}_2$. Hence
${\tilde W} (k, p; p_{0,2})$ provides analytic continuation of
${\hat W} (k, p)$ to the interval $\mathcal{I}_2$. We can continue
this process to obtain analytic continuation of ${\hat W}$ to any
interval $\mathcal{I}_j$. Since the union of $\left \{ \mathcal{I}_j \right \}_{j=0}^\infty$
contains $\mathbb{R}^+ \cup \{0 \}$, it follows that
${\hat W} (k, .)$ is analytic in $\mathbb{R}^+$.  
In particular, (\ref{100}) provides the local
Taylor series representation of ${\hat W} (k, p)$ near $p=p_0$.
\end{proof}

\noindent{\bf Proof of Theorem \ref{T2}}  

Using Lemma \ref{lem40}, it follows from
the inequality
$ \|W^{[l]} (\cdot, p_0) \|_{\infty} \le 
\| {\hat W}^{[l]} (\cdot, p_0) \|_{\mathbb{L}^1} $ by
integration in $k$ that 
$$ |W^{[l]} (x, p_0)| \le \frac{8 \pi A (4 B)^l e^{\alpha p_0}}{
\beta (2l+1)^2 (1+p_0^2)} $$
$$ |D W^{[l]} (x, p_0)| \le \frac{8 \pi A (4 B)^l e^{\alpha p_0}}{
\beta (2l+1)^2 (1+p_0^2) } $$
$$ |D^2 W^{[l]} (x, p_0)| \le \frac{16 \pi A (4 B)^l e^{\alpha p_0} }{
\beta^2 (2l+1)^2 (1+p_0^2)}
$$
and therefore, the series (\ref{100})
converges for $|p-p_0| < B^{-1}/4$ and, from Lemma 
\ref{lemancont} it is the 
local representation of the solution ${\hat W} (k, p)$ to 
(\ref{100.3}) satisfying ${\hat W} (k, 0) =0$
for any $p_0 \ge 0$. These estimates  
on $W $ in terms of ${\hat W}$,
and the fact that $W (x, p) $ is analytic in a neighborhood of
for $p \in \left \{ 0 \right \} \cup \mathbb{R}^+$ and is
exponentially bounded in $p$ for large $p$ 
(recall ${\hat W} \in \mathcal{A}^\alpha$) implies
Borel summability of $v$ in  $1/t$.
Watson's Lemma \cite{Wasow} implies
$w (x, t) = \int_0^\infty e^{-p/t} W(x, p) dp 
\sim \sum_{m=2}^\infty v_m (x) t^m $,
implying
$$ v(x, t) = v_0 (x) + t v_1 (x) + 
\sum_{m=2}^\infty v_m (x) t^m, 
$$
where $ v_m (x) = m! W^{[m-1]} (x; 0) = m! U^{[m-1]} (x; 0)$ for $m\ge 2$.
It follows from the bounds on
${\hat W}^{[m-1]} (k)$ in \S 3, that for $m \ge 2$,
$ | W^{[m-1]} (x; 0) | \le A_0 B_0^m$, where $A_0$ and $B_0$\footnote{
We may express
it in terms of $A$ and $B$ as well, however, the estimates $A_0$ and
$B_0$  found in \S 3, are better.}
are chosen 
to ensure (\ref{12.4.0}), (\ref{ab2cond}) and
(\ref{Bcond}). 

\section{Acknowledgments.} 
The authors benefitted from comments by Peter Constantin and Charlie Doering.
This work was supported in part by the National
Science Foundation (DMS-0406193, DMS-0601226, DMS-0600369 to OC and
(DMS-0405837 to S.T). Additional support was provided to ST by the Institute
for Mathematical Sciences, Imperial College and the EPSRC.

\section{Appendix}

\subsection{Some Fourier convolution inequalities}
The following lemmas are relatively straightforward.

\begin{Definition}  
\label{defa1}
Consider the polynomial
$$ P_n (z) = 
\sum_{j=0}^n \frac{n!}{j!} z^{j} $$  
\end{Definition}
\begin{Remark}
Integration by parts yields
\begin{equation}
\label{eqintPn}
\int_{0}^z e^{-\tau} \tau^n d\tau = - e^{-z} P_n (z) + n!  
\end{equation}
\end{Remark}
\begin{Lemma}
\label{lema2}
For all $y \ge 0$ and nonnegative integers $m, n \ge 0$ we have
$$ y^{m+1} \int_0^1 \rho^{m} P_n (y (1-\rho)) d\rho  =
m! n! \sum_{j=0}^n \frac{y^{m+j+1}}{(m+j+1)!} 
$$
\end{Lemma}

\begin{proof}
This follows from a simple computation:
$$y^{m+1} \int_0^1 \rho^m P_n (y (1-\rho))d\rho = 
\sum_{j=0}^n \frac{n!}{j!} y^{j+m+1} \int_0^1 (1-\rho)^{j} \rho^m d\rho 
= 
m! n! \sum_{j=0}^n \frac{y^{j+m+1}}{(m+j+1)!}  
$$
\end{proof}

\begin{Lemma}
\label{lema1.1}
For all $y \ge 0$ and nonnegative integers $n\ge m \ge 0$ we have
$$ y^{m+1} \int_1^\infty e^{-2 y(\rho-1)} \rho^m P_n (y (\rho-1)) d\rho 
\le 2^{-m} (m+n)! \sum_{j=0}^m \frac{y^j}{j!}   
$$
\end{Lemma}
\begin{proof}
First we note that
\begin{multline*}
y^{m+1+l} \int_1^\infty e^{-2 y(\rho-1)} \rho^m (\rho-1)^l d\rho
= 
y^{m+1+l} \int_0^\infty e^{-2 y \rho} (1+\rho)^m \rho^l d\rho \\
= y^{m+1+l} \sum_{j=0}^m \frac{m!}{j! (m-j)!} \int_0^\infty e^{-2 y \rho}
\rho^{l+j} d\rho =  
2^{-l-1} \sum_{j=0}^m \frac{y^{m-j} m! (l+j)!}{j! (m-j)! 2^{j}} \\
= 2^{-l-1} \sum_{j=0}^m \frac{y^{j} m! (l+m-j)!}{(m-j)! j! 2^{m-j}} 
\end{multline*}
Therefore, from the definition of $P_n$, it follows that
\begin{multline*}
y^{m+1} \int_1^\infty e^{-2 y(\rho-1)} \rho^m P_n (y (\rho-1)) d\rho 
\\= m! n! \sum_{j=0}^m \frac{y^j}{j! (m-j)! 2^{m-j}} 
\left ( \sum_{l=0}^n  
\frac{(l+m-j)!}{2^{l+1} l!} \right ) 
\end{multline*}
Taking the ratio of two consecutive terms we see that ${(l+m-j)!}/{l!}$ is
nondecreasing with $l$ since $m-j \ge 0$. Therefore the $l=n$ term is the
largest term in the summation over $l$. Further, $\sum_{l=0}^n 2^{-l-1} \le
1$. Therefore, $ \sum_{l=0}^n 2^{-l-1}{(l+m-j)!}/{l!} \le
{(m-j+n)!}/{n!} $, and hence
$$
y^{m+1} \int_1^\infty e^{-2 y(\rho-1)} \rho^m P_n (y (\rho-1)) d\rho 
\le 
2^{-m} m! n! \sum_{j=0}^m \frac{y^j}{j!} \frac{2^j (m-j+n)!}{n! (m-j)!} $$
The ratio of two consecutive (in $j$) terms in ${2^j
  (m-j+n)!}/{(m-j)!}$ is $\le 1$ for  $m \le n$, hence
the largest value is attained at $j=0$ and thus
$$y^{m+1} \int_1^\infty e^{-2 y(\rho-1)} \rho^m P_n (y (\rho-1)) d\rho 
\le 
2^{-m} (m+n)! \sum_{j=0}^m \frac{y^j}{j!} $$
\end{proof}

\begin{Lemma}
\label{lema1.2}For all $y \ge 0$ and nonnegative integers $n\ge m \ge 0$
we have
$$ y^{m+1} \int_0^\infty e^{-y (\rho-1) [1+\text{sgn} (\rho-1)]}
\rho^{m} P_n (y |1-\rho|) d\rho  \le 
m! n! Q_{m+n+1} (y) 
$$
\end{Lemma}
\begin{proof}
By breaking up the integral range into $\int_0^1$ and
$\int_{1}^\infty$ and using the two previous Lemmas, we
obtain
\begin{multline*}
y^{m+1} \int_0^\infty e^{-y (\rho-1) [1+\text{sgn} (\rho-1)]}
\rho^{m} P_n (y |1-\rho |) d\rho  \le 
m! n! \Bigg ( \sum_{j=m+1}^{m+n+1} \frac{y^j}{j!} \\+ 
2^{-m-n} \frac{(m+n)!}{m! n!} \sum_{j=0}^m 2^{n} 
\frac{y^j}{j!} \Bigg) 
\le 
m! n! \sum_{j=0}^{m+n+1} 2^{m+n+1-j} \frac{y^j}{j!} = m! n! Q_{m+n+1} (y) 
\end{multline*}
where we used
$2^{-m-n} \frac{(m+n)!}{m! n!} \le 1$.   
\end{proof}

\begin{Lemma}
\label{lema0.0}
If $m$ and $ n $, are integers no less that $-1$ we obtain
$$
|q| \int_{q'\in \mathbb{R}^3} 
e^{|q|- |q'| - |q-q'|} |q'|^{m} |q-q'|^{n} dq'     
\le 2 \pi (m+1)! (n+1)! Q_{m+n+3} (|q|)
$$     
\end{Lemma}

\begin{proof}
We note that we may assume $m \le n$ without loss of generality
since   changing variable
$q'\mapsto q-q'$  switches the roles of $m$ and $n$. 

First, we will show that
\begin{multline}
\label{eqshow}
\frac{|q|}{2\pi} \int_{q'\in \mathbb{R}^3} 
e^{|q|- |q'| - |q-q'|} |q'|^{m} |q-q'|^{n} dq'     
\\\le |q|^{m + 2} \int_0^\infty e^{- |q| (\rho-1) [1+{\text{sgn}}(\rho-1)]}
\rho^{m+1} P_{n+1} (|q| (|\rho-1|) d\rho  
\end{multline}
We scale $q'$ with $|q|$ and use a polar representation $(\rho, \theta, \phi)$
for $q'/|q|$, where $\theta$ is the angle between $q$ and $q'$. As a variable
of integration however, we prefer to use $z =
\sqrt{1+\rho^2-2 \rho \cos \theta}$ to $\theta$.
Then, it
is clear that 
$$|q-q'| = |q| \sqrt{1+\rho^2 - 2 \rho \cos \theta} = |q|
z, ~~~{\rm and} ~~~~~dz = \frac{\rho \sin \theta d\theta}{
\sqrt{1+\rho^2-2 \rho \cos \theta}} $$ 
Therefore,  
\begin{multline*}
|q| \int_{q' \in \mathbb{R}^3} 
e^{|q|- |q'| - |q-q'|} |q'|^{m} |q-q'|^{n} dq' 
\\
= 2 \pi |q|^{m + n + 4} \int_0^\infty d\rho \rho^{m+1} e^{-|q| (\rho-1)} 
\left \{
\int_{|\rho-1|}^{1+\rho} dz e^{-|q| z} z^{n+1} \right \} \\          
= 2 \pi |q|^{m+2} \int_0^\infty d\rho \rho^{m+1} e^{-|q| (\rho-1)}
\left [ e^{-|q| |\rho-1|}
P_{n+1} (|q| |\rho-1|)  - e^{-|q| (1+\rho)} P_{n+1} (|q| (1+\rho)) \right ] 
\end{multline*}
Inequality (\ref{eqshow}) follows since
$e^{-|q| (1+\rho)} P_{n+1} (|q| (1+\rho)) \ge  0$.
The rest of the Lemma follows from 
Lemma \ref{lema1.2}, with  $y=|q|$, and  $m$ replaced  by $m+1$, 
$n$ by $n+1$ respectively.
\end{proof}

\begin{Lemma}
\label{lema1.3}
For any $\mu \ge  1 $, and nonnegative integers $m, n$ we have
\begin{multline*}
|k| \int_{k'\in \mathbb{R}^3} 
\frac{e^{- \beta [|k'| + |k-k'|] } }{ (1+|k'|)^\mu (1+|k-k'|)^\mu} 
|\beta k'|^{m} |\beta (k-k')|^{n} dk' 
\\
\le 
\frac{\pi 2^{\mu+1} e^{-\beta |k|} m! n!}{\beta^3 (1+|k|)^\mu} 
( m+n+2 )
Q_{m+n+2} (\beta |k|)       
\end{multline*}
\end{Lemma}
\begin{proof}
We break up the integral into two ranges: 
\begin{equation}
\label{eqbreak}
\int_{|k'| \le |k|/2} + \int_{|k'| > |k|/2}
\end{equation} 
In the first integral we have 
$$\frac{1}{(1+|k-k'|)^\mu (1+|k'|)^\mu } 
\le \frac{1}{(1+|k|/2)^\mu (1+|k'|)} \le    
\frac{\beta }{(1+|k|/2)^\mu |\beta k'|}$$
While in the second integral we have
$$\frac{1}{(1+|k-k'|)^\mu (1+|k'|)^\mu } 
\le \frac{1}{(1+|k|/2)^\mu (1+|k-k'|)} \le    
\frac{\beta }{(1+|k|/2)^\mu |\beta (k-k')|}$$
Introducing in the first integral
$q= \beta k$ and $q' = \beta k'$, we obtain  
\begin{multline*}
|k| \int_{|k'| \le |k|/2}  
\frac{ e^{- \beta [|k'| + |k-k'|]} }{(1+|k'|)^\mu (1+|k-k'|)^\mu} 
|\beta k'|^{m} |\beta (k-k')|^{n} dk' \\
\le  
\frac{2^\mu e^{-\beta |k|} }{\beta^3 (1+|k|)^\mu}     
|q| \int_{q' \in \mathbb{R}^3} 
e^{- |q'|- |q-q'| + |q|}
|q'|^{m-1} |q-q'||^{n} dq' 
\end{multline*}
while in the second integral, with
$q=\beta k$ and  $q-q'=\beta k'$, we obtain
\begin{multline*}
|k| \int_{|k'| > |k|/2}  
\frac{ e^{- \beta [|k'| + |k-k'|]} }{(1+|k'|)^\mu (1+|k-k'|)^\mu} 
|\beta k'|^{m} |\beta (k-k')|^{n} dk' \\ 
\le  
\frac{2^\mu  e^{-\beta |k|} }{ \beta^3 (1+|k|)^\mu}     
|q| \int_{q' \in \mathbb{R}^3} 
e^{- |q'|- |q-q'| + |q|}
|q'|^{n-1} |q-q'|^{m} dq' 
\end{multline*}
We now use 
Lemma \ref{lema0.0} to bound the first integral, with $m$ 
replaced by 
$m-1$.  We also use
Lemma \ref{lema0.0} to bound the second integral, 
with $n-1$ replacing $n$.
The proof is completed by adding the two bounds.
\end{proof}

\begin{Lemma}
\label{lema1.3.0}
For any $\mu \ge  2 $, and  $n \in\NN\setminus\{0\} $ we have
\begin{multline*}
|k| \int_{k'\in \mathbb{R}^3} 
\frac{e^{- \beta [|k'| + |k-k'|] } }{ (1+|k'|)^\mu (1+|k-k'|)^\mu} 
|\beta (k-k')|^{n} dk' 
\\
\le 
\frac{2^{\mu+1} \pi e^{-\beta |k|} }{\beta^2 (1+|k|)^\mu}  
\left \{ (n-1)! Q_{n+1} (|q|) + \frac{3 (n+1)! |q|^{2/3}}{2 \beta^{2/3}} 
\sum_{j=0}^{n+1} \frac{|q|^j}{j!} \right \}  
\end{multline*}
\end{Lemma}
\begin{proof}
We break up the integral into $\int_{|k'| < |k|/2} + \int_{|k'| \ge |k|/2}$. 
In the first integration range we have
$[1 + |k-k']^{-\mu} \le 2^\mu [1+|k|]^{-\mu}$, whereas in the second
range $[1+|k'|]^{-\mu} \le 2^\mu [1+|k|]^{-\mu}$.  Therefore,
using Lemma \ref{lema0.0},
\begin{multline*}
|k| \int_{|k'| \ge |k|/2}  
\frac{e^{- \beta [|k'| + |k-k'|] } }{ (1+|k'|)^\mu (1+|k-k'|)^\mu} 
|\beta (k-k')|^{n} dk' 
\\
\le 
\frac{2^{\mu} e^{-\beta |k|} }{\beta^2 (1+|k|)^\mu} 
|q| \int_{q \in \mathbb{R}^3} e^{-|q'|-|q-q'|+|q|} ~|q-q'|^{n-2} dq'  
\le 
\frac{2^{\mu+1} \pi e^{-\beta |k|} }{\beta^2 (1+|k|)^\mu}  
(n-1)! 
Q_{n+1} (|q|) 
\end{multline*}
On the other hand, using $(1+|k'|)^{-\mu} \le (1+|k'|)^{-2+2/3} \le 
|k'|^{-2+2/3}$ we get 
\begin{multline*}
|k| \int_{|k'| < k|/2}  
\frac{e^{- \beta [|k'| + |k-k'|] } }{ (1+|k'|)^\mu (1+|k-k'|)^\mu} 
|\beta (k-k')|^{n} dk' \\
\le  
\frac{2^\mu e^{-\beta |k|} }{\beta^{2+2/3} (1+|k|)^\mu} 
|q| \int_{|q'| \le |q|/2}
|q'|^{-2+2/3} 
e^{-|q'|-|q-q'|+|q|} ~|q-q'|^{n} dq'  
\end{multline*}
We note that 
\begin{multline*}
\frac{|q|}{2\pi} \int_{|q'| < |q|/2}
 |q'|^{-2+2/3} 
e^{-|q'|-|q-q'|+|q|} ~|q-q'|^{n} dq' \\ 
= |q|^{n+2+2/3} \int_0^{1/2} \rho^{-1+2/3} e^{-|q| (\rho - 1)} 
\left \{
\int_{1-\rho}^{1+\rho} dz e^{-|q| z} z^{n+1} 
\right \} d\rho \\
\le |q|^{2/3} \int_0^{1/2} \rho^{-1+2/3} P_{n+1} (|q| (1-\rho) 
d\rho 
\\ \le |q|^{2/3} \sum_{j=0}^{n+1} \frac{(n+1)!}{j!} |q|^j \int_0^{1} 
\rho^{-1/3} (1-\rho)^j d\rho  
\le 
\frac{3}{2} |q|^{2/3} (n+1)! \sum_{j=0}^{n+1} 
\frac{|q|^j}{j!}
\end{multline*}
\end{proof}
\begin{Lemma}
\label{lema1.4}
For any $\mu \ge 1$ and nonnegative integers $l_1, l_2 \ge 0$ we have 
\begin{multline}
\label{lem6.9eq}
|k| \int_{k'\in \mathbb{R}^3} 
\frac{ e^{- \beta [|k'| + |k-k'|] } }{ (1+|k'|)^\mu (1+|k-k'|)^\mu} 
Q_{2l_1} (|\beta k'|) Q_{2l_2} (|\beta (k-k')|) dk' \\
\le \frac{2^{\mu+1}\pi e^{-\beta |k} }{ 3 \beta^3 (1+|k|)^\mu}
(2l_1+2l_2+1) (2l_1+2l_2+2) (2l_1+2l_2+3) Q_{2l_1+2l_2+2} (\beta |k|)   
\end{multline}
\end{Lemma}
\begin{proof}
As before, we define $q=\beta k$.
Also, for notational convenience, we define 
$$ k^m \Join k^n 
=   
|k| \int_{k'\in \mathbb{R}^3} 
\frac{ e^{- \beta [|k'| + |k-k'|] } }{ (1+|k'|)^\mu (1+|k-k'|)^\mu} 
|\beta k'|^m |\beta (k-k')|^n dk' 
$$
$$ K = \frac{2^{\mu+1}\pi e^{-\beta |k|}}{\beta^3 (1+|k|)^\mu} $$ 
Lemma \ref{lema1.3} and 
$2^{2l_1+2l_2+2-(j+2)} Q_{j+2} (|q|) \le Q_{2l_1+2l_2+2}$ for 
$j \le 2 l_1 + 2l_2$ imply 
that the left side of (\ref{lem6.9eq}) is given by
\begin{multline*}
\sum_{m=0}^{2l_1} \sum_{n=0}^{2l_2}  
\frac{2^{2l_1+2l_2-m-n}}{m! n!} k^m \Join k^n  \le 
K \sum_{m=0}^{2l_1} \sum_{n=0}^{2l_2}  
2^{2l_1+2l_2-m-n} (m+n+2)  
Q_{m+n+2} (\beta |k|) \\
\le 
K \sum_{j=0}^{2l_1+2l_2} 2^{2l_1+2l_2+2-(j+2)} 
(j+2) (j+1) Q_{j+2} (|q|) \\
\le K Q_{2l_1+2l_2+2} (|q|) 
\sum_{j=0}^{2l_2+2l_1} (j+1) (j+2) \\\le    
\frac{K}{3} (2l_1+2l_2+1)(2l_1+2l_2+2)(2l_1+2l_2+3) Q_{2l_1+2l_2+2} (|q|), 
\end{multline*}  
which imply the result.
\end{proof}
\begin{Lemma}
\label{lemaQ2l}
If $ \mu \ge 2$ and $l \ge 0$, 
then 
\begin{multline*}
\frac{|k|}{(l+1)^{2/3}} 
\int_{k'\in \mathbb{R}^3} 
\frac{ e^{- \beta [|k'| + |k-k'|] } }{ (1+|k'|)^\mu (1+|k-k'|)^\mu} 
Q_{2l} (|\beta (k-k')|) dk' \\
\le \frac{C_1 e^{-\beta |k} }{
(1+|k|)^\mu} (2l+1) Q_{2l+2} (\beta |k|),
\end{multline*}
where
\begin{equation*}
C_1 = 
12 \pi 2^\mu \beta^{-8/3} + 2 \pi 2^\mu \beta^{-2} + \frac{1}{2} C_0 (\mu) 
\beta^{-1} 
\end{equation*}
\end{Lemma}
\begin{proof}
The case $l=0$ follows easily by using (\ref{eqlem0.1}) and the fact
that 
$$|k| = \beta^{-1} |q| \le \frac{1}{2} \beta^{-1} Q_{2} (|q|) $$ 
For $l \ge 1$,
it is convenient to separate out the constant term $2^{2l}$ in $Q_{2l}$ 
and note that from (\ref{eqlem0.1}) and the definition of $Q_{n} (z)$ we have
\begin{equation*}
|k| \int_{k'\in \mathbb{R}^3} 
\frac{e^{-\beta [|k'| + |k-k'|] } }{ (1+|k'|)^\mu (1+|k-k'|)^\mu} 
2^{2l} dk' \le 
\frac{C_0 |k| e^{-\beta |k|} }{(1+|k|)^\mu } 
2^{2l} 
\le 
\frac{C_0 e^{-\beta |k|} }{ 2 \beta (1+|k|)^\mu } 
Q_{2l} (\beta |k|) 
\end{equation*}  
As in previous Lemma, for notational convenience, we define 
$$ k^m \Join k^n 
=   
|k| \int_{k'\in \mathbb{R}^3} 
\frac{ e^{- \beta [|k'| + |k-k'|] } }{ (1+|k'|)^\mu (1+|k-k'|)^\mu} 
|\beta k'|^m |\beta (k-k')|^n dk' 
$$
Then, it is clear from Lemma \ref{lema1.3.0} that 
\begin{multline*}
[Q_{2l} (\beta |k|)-2^{2l} ] \Join k^0 =
\sum_{n=1}^{2l} \frac{2^{2l-n}}{n!} k^n \hat{*} k^0 \\
\le    
\frac{2^{\mu+1} \pi e^{-\beta |k|}}{\beta^2 (1+|k|)^\mu} 
\sum_{n=1}^{2l} \frac{2^{2l-n}}{n!} 
\left \{ (n-1)! \sum_{j=0}^{n+1} \frac{2^{n+1-j} (\beta |k|)^j}{j!}  
+ \frac{3 (n+1)! |q|^{2/3}}{2 \beta^{2/3}} 
\sum_{j=0}^{n+1} \frac{(\beta |k|)^j}{j!} \right \}  \\
\le \frac{2^{\mu+1} \pi e^{-\beta |k|}}{\beta^2 (1+|k|)^\mu} \Bigg \{
\sum_{j=0}^{2l+1} \frac{2^{2l+1-j} (\beta |k|)^j}{j!} 
\sum_{n'=\max\{j,2\} }^{2l+1} \frac{(n'-2)!}{(n'-1)!}     
\\+\frac{3 |q|^{2/3}}{\beta^{2/3} 2}  
\sum_{j=0}^{2l+1} \frac{2^{2l+1-j} (\beta |k|)^j}{j!} 
\sum_{n'=\max\{j,2\}}^{2l+1} 2^{j-n'} n' \Bigg \} \\
\le \frac{2^{\mu+1} \pi e^{-\beta |k|}}{\beta^3 (1+|k|)^\mu} \left [
\beta Q_{2l+1} (\beta |k|) \log (2l+2) 
+ 3 (2l+1)  
\beta^{1/3} |\beta k|^{2/3} Q_{2l+1} (\beta |k|)
\right ] 
\end{multline*}
The lemma follows since $\log (2l+2) / (2l+1) \le 1 $, while
if $ |\beta k| \le (l+1)$,  
$$ \left ( \frac{|\beta k|}{(l+1)} \right )^{2/3} Q_{2l+1} (\beta |k|) 
\le Q_{2 l+1} (\beta |k|) \le \frac{1}{2} Q_{2l+2} (\beta |k|)$$ 
whereas for $\beta |k| \ge (l+1)$ we have   
$$ \left ( \frac{|\beta k|}{(l+1)} \right )^{2/3} Q_{2l+1} (\beta |k|) 
\le \frac{2 \beta |k|}{2 l+2} Q_{2l+1} (\beta |k|) 
\le 2 Q_{2l+2} (\beta |k|)$$
\end{proof}

\vfill \eject
\end{document}